%% file: ChaikovskiiZhang.tex
\documentclass[onefignum,onetabnum]{siamonline220329}

\input{ex_shared}

\ifpdf
\hypersetup{
  pdftitle={Solving forward and inverse problems in a non-linear 3D PDE via an asymptotic expansion based approach},
  pdfauthor={Dmitrii Chaikovskii and Ye Zhang}
}
\fi

\begin{document}

\maketitle

\begin{abstract}
This paper concerns the use of asymptotic expansions for the efficient solving of forward and inverse problems involving a nonlinear singularly perturbed time-dependent reaction--diffusion--advection equation. By using an asymptotic expansion with the local coordinates in the transition-layer region, we prove the existence and uniqueness of a smooth solution with a sharp transition layer for a three-dimensional partial differential equation. Moreover, with the help of asymptotic expansion, a simplified model is derived for the corresponding inverse source problem, which is close to the original inverse problem over the entire region except for a narrow transition layer. We show that such simplification does not reduce the accuracy of the inversion results when the measurement data contain noise. Based on this simpler inversion model, an asymptotic-expansion regularization algorithm is proposed for efficiently solving the inverse source problem in the three-dimensional case. A model problem shows the feasibility of the proposed numerical approach.
\end{abstract}



\section{Introduction} \label{sec:Introduction}

The main goal of this paper is to develop a new framework for solving the following three-dimensional inverse source problem efficiently.

(\textbf{IP}): Given full noisy data $\{u^\delta, u_x^\delta, u_y^\delta, u_z^\delta\}$ of $\{u,u_x,u_y,u_z\}$ or partial noisy data $\{u^\delta\}$ of $\{u\}$ at the $p\cdot q \cdot v$ location points $\{x_i, y_j, z_k\}^{p,q,v}_{i,j,k=0}$ and at the time point $t_0$, find the source function $f(x,y,z)$ such that $(u,f)$ satisfies the following dimensionless nonlinear problem:
\begin{align} \label{mainproblem}
\begin{dcases}
\displaystyle \mu \Delta u-\frac{\partial u}{\partial t} =-u \left( \frac{\partial u}{\partial x}+ \frac{\partial u}{\partial y} + \frac{\partial u}{\partial z} \right)+f, \ (x, y) \in \mathbb{R}^2, \ z \in (-a, a) \equiv \Omega , \ t\in (0, T] \equiv \mathcal{T}, \\
u(x, y,z, t, \mu)=u(x+L, y,z, t, \mu)=u(x, y+M,z, t, \mu), \quad(x, y,z,t) \in \mathbb{R}^2 \times \bar{\Omega} \times \bar{\mathcal{T}}, \\
u(x, y,-a, t, \mu)=u^{-a}(x,y), \quad u(x, y,a, t, \mu)=u^{a}(x,y), \quad(x, y ,t) \in \mathbb{R}^2 \times \bar{\mathcal{T}}, \\
u(x, y,z, 0, \mu)=u_{init}(x, y,z, \mu), \quad(x, y,z) \in \mathbb{R}^2 \times \bar{\Omega},
\end{dcases}
\end{align}
where $u(x,y,z,t)$ represents the dimensionless value (e.g., of the temperature), $\mu \ll 1$ is a small parameter (also called the diffusion coefficient), and $f(x,y,z)$ is the source function. Suppose that (i) the function $f(x,y,z)$ is $L$-periodic in the variable $x$, $M$-periodic in the variable $y$, and sufficiently smooth in the region $(x,y,z): \mathbb{R}^2 \times \bar{\Omega}$ ($\Omega\equiv(-a, a)$), (ii) the functions $u^{-a}(x,y)$ and $u^{a}(x,y)$ are $L$-periodic in $x$, $M$-periodic in $y$, and sufficiently smooth in $(x,y)\in \mathbb{R}^2$, and (iii) $u_{init}(x,y,z,\mu)$ is a sufficiently smooth function in $(x,y,z):\mathbb{R}^2 \times \Omega$, is $L$-periodic in $x$ and $M$-periodic in $y$, and satisfies $u_{init}(x,y,-a,\mu)=u^{-a}(x,y), u_{init}(x,y,a, \mu)=u^{a}(x,y)$.

The reaction--diffusion--advection (RDA) system with small parameters plays an important role in the quantitative study of many scientific problems, including chromatography \cite{zhang2016regularization,LinZhang2018}, biology \cite{HidalgoTello2014}, and environmental problems \cite{2014evolutiondisp}. This paper is focused on the solutions of RDA system \eqref{mainproblem} with a large gradient in a certain region. This region is called the inner transition layer, whose width is usually small compared to that of the entire region. The presence of a small parameter accompanying the highest derivative makes the equation singularly perturbed, i.e., one cannot simply ignore the term with a small parameter in the differential equation because doing so would lead to a totally different solution that would not reflect the physical phenomenon being studied. RDA problems with inner transition layers are often used for the mathematical modeling of the density distributions of liquids or gases or temperature in the presence of spatial inhomogeneities \cite{Sogachev2006ModificationOT, Olchev2009ApplicationOA} or in nonlinear acoustics \cite{Rudenko2017IBE, Nefedov2018OFM, Nefedov2019TheEA}, and such problems include the temperature behavior in the near-surface layer of the ocean \cite{LevNik17}, the carrier wave functions in Si/Ge heterostructures \cite{Orlov_2015}, and propagation of an autowave front in a medium with barriers \cite{b04, Sidorova2019AutowaveMO}. Determining existence conditions for and the stability of nonstationary solutions with large gradients is important for creating adequate models of processes with nonstationary distributions of fields of physical quantities \cite{Kopteva2011StabilisedAO, VoNEAN2015}. Analytical studies also make it possible to create efficient numerical methods for solving equations with inner transition layers \cite{IJNAM-7-393, QUINN2015500, ORiordan2011ParameteruniformNM, Lukyanenko2017, LukyanenkoTwoDimInv2018}. The technique of asymptotic analysis is remarkably good for overcoming the aforementioned difficulties, i.e., the issues of (i) existence and uniqueness of smoothing solutions and (ii) high qualified resolution of approximate analytical solutions.

On the other hand, because the inverse source problem (\textbf{IP}) is ill-posed (i.e., small noise in the measurement data may lead to arbitrarily large changes in the classical approximate solution, which is far from the ground truth; see \cite{TikhonovYagola1995,Isakov1990} for details), a regularization technique should be developed for a stable solution of (\textbf{IP}). Following the framework of Tikhonov regularization, (\textbf{IP}) can be converted to the following minimization problem with partial differential equation (PDE) constraints:
\begin{multline}
\label{LSM0}
\min_{f} \sum^n_{i=0} \sum^m_{j=0} \left\{ \left[u(x_i,y_j,z_k,t_0) - u^\delta(x_i,y_j,z_k,t_0) \right]^2 \vphantom{\left[\frac{\partial u}{\partial x} \right]^2 }\right.\\
+ \left[\frac{\partial u}{\partial x}(x_i,y_j,z_k,t_0) - \frac{\partial u^\delta}{\partial x}(x_i,y_j,z_k,t_0) \right]^2 + \left[\frac{\partial u}{\partial y}(x_i,y_j,z_k,t_0) - \frac{\partial u^\delta}{\partial y}(x_i,y_j,z_k,t_0) \right]^2\\
\left. + \left[\frac{\partial u}{\partial z}(x_i,y_j,z_k,t_0) - \frac{\partial u^\delta}{\partial z}(x_i,y_j,z_k,t_0) \right]^2 \right\} + \varepsilon \mathcal{R}(f),
\end{multline}
where $u$ solves the nonlinear PDE \eqref{mainproblem} with a given $f$, $\mathcal{R}(f)$ is the regularization term, and $\varepsilon>0$ is the regularization parameter.

There are two essential difficulties when employing the solution model \eqref{LSM0}. First, it is difficult to select an appropriate regularization term $\mathcal{R}$ and regularization parameter $\varepsilon$ in practice. By the standard argument of the regularization theory of inverse problems, the optimal choices of these two quantities depend on the ground truth $f$, which is unknown in real-world problems. Second, even if both $\mathcal{R}$ and $\varepsilon$ are given, the numerical realization of the PDE-constrained optimization problem \eqref{LSM0} is a very hard task because it is a nonconvex optimization with a nonlinear PDE constraint. Therefore, the main purpose of this work is to find a new model to replace \eqref{LSM0}, one that is much simpler but still sufficiently accurate. Fortunately, this can be achieved by using the theory of asymptotic analysis, which is also a main motivation of this work.

Indeed, several previous studies have used asymptotic analysis to solve inverse problems involving singularly perturbed PDEs; for instance, in \cite{Lukyanenko2018SolvingOT}, a coefficient inverse problem for the one-dimensional RDA equation was solved using the position of the moving front, which was obtained approximately using asymptotic analysis, and in \cite{LukyanenkoTwoDimInv2018}, an inverse problem of determining the position of a moving front for a two-dimensional reaction--diffusion equation was studied via asympototic analysis. Actually, the present paper can be viewed as an extension of our previous work \cite{CHAIKOVSKII2022111609}, which dealt with one-dimensional inverse source problems under some a~priori structure assumptions about the source function. The novelty of the present work is combining asymptotic expansions and local coordinates to solve the forward problem, which in turn allows us to reduce the inverse problem of finding the source function in the original PDE with a large derivative to a simpler equation and obtain sufficiently accurate results. Although local coordinates have been used previously to study some two-dimensional forward \cite{AntLevNef18} and inverse \cite{LukyanenkoTwoDimInv2018} problems, to the best of our knowledge this is the first time that the combination of local coordinates and Vasil'eva's algorithm for asymptotic expansion in a small parameter \cite{Vasil1990} has been used. Moreover, the mathematical analysis in a three-dimensional setting---including the existence and uniqueness of the forward problem and the inversion method for (\textbf{IP})---has never been investigated before.

This paper is organized as follows. Section~\ref{constructionOfAsymptotic} presents the asymptotic analysis for the forward problem. Section~\ref{statementresults} presents the main theoretical results, while Section~\ref{derivationAndProofs} presents their technical proofs. Section~\ref{simulation} presents some experiments for the forward and inverse problems. Finally, Section~\ref{Conclusion} presents concluding remarks.

\section{Asymptotic analysis} \label{constructionOfAsymptotic}

We investigate the solution of problem \eqref{mainproblem}, which has the form of a moving front that at each moment of time is close to $\varphi^{(-)}(x,y,z)$ for $-a\leq z\leq h(x,y,t)$, close to $\varphi^{(+)}(x,y,z)$ for $h(x,y,t)\leq z\leq a$, and changes sharply from $\varphi^{(-)}(x,y,z)$ to $\varphi^{(+)}(x,y,z)$ in a neighbourhood of the surface $z=h(x,y,t)$. In this case, the solution to problem \eqref{mainproblem} has an inner transition layer in the vicinity of this surface.

First, we list all the main assumptions used throughout the paper.

\begin{assumption}
\label{A1}
$u^{-a}(x,y)<0$, $u^{a}(x,y)>0$, and $u^{a}(x,y)-u^{-a}(x,y)> 2\mu^2 $ for all $(x,y) \in \mathbb{R}^2$.
\end{assumption}

\begin{assumption}
\label{A2}
For any $(x,y,z) \in \mathbb{R}^2\times \bar{\Omega}$,
\begin{align*}
\left(u^{-a}\left({\scriptstyle x-a-z, y-a-z }\right)\right)^2 &> -2 \int _0^{a+z} f\left({\scriptstyle x-a-z+s, y-a-z+s, -a+s }\right)ds, \\
\left(u^{a}\left({\scriptstyle x+ a-z, y+a-z }\right)\right)^2 &> 2 \int _{z-a}^{0} f\left({\scriptstyle x+a-z+s, y+a-z+s, a+s }\right)ds.
\end{align*}
\end{assumption}

\begin{assumption}\label{A3}
$-a < {h_0}(x,y, t)<a$ for any $(x,y,t) \in \mathbb{R}^2\times \mathcal{T}$, where $h_0$ is the zero approximation of $h$ [see \eqref{curveexpansion}] and $ \displaystyle \max_{(x,y)\in \mathbb{R}^2,\ t\in \bar{\mathcal{T}} } \left( \frac{\partial {h_0} (x,y,t)}{\partial x}+ \frac{\partial {h_0} (x,y,t)}{\partial y} \right) <1 $.
\end{assumption}

\begin{assumption}
\label{A4}
$u_{init}(x,y,z)=U_{n-1}(x,y,z,0)+\mathcal{O}(\mu^n)$, where the asymptotic solution $U_{n-1}$ is constructed later; see Theorem~\ref{MainThm}, for example.
\end{assumption}

Under Assumptions~\ref{A1} and \ref{A2}, $\varphi^{(-)}(x,y,z)$ and $\varphi^{(+)}(x,y,z)$---which are used to construct the zeroth-order asymptotic solution [cf.~\eqref{u0regu}]---can be calculated explicitly from the reduced stationary equations \eqref{zeroorderregularequation1} and \eqref{zeroorderregularequation2} by using the method of characteristics:
\begin{align}
\label{eq7}
\displaystyle \begin{dcases}
\varphi^{(-)} &= -\sqrt{\left(u^{-a}\left({\scriptstyle x-a-z, y-a-z} \right)\right)^2+2 \int _0^{a+z} f\left( {\scriptstyle x-a-z+s, y-a-z+s, -a+s }\right)ds},\\
\varphi^{(+)} &= \sqrt{\left(u^{a}\left({\scriptstyle x+ a-z, y+a-z }\right)\right)^2-2 \int _{z-a}^{0} f\left({\scriptstyle x+a-z+s, y+a-z+s, a+s }\right)ds}.
\end{dcases}
\end{align}

The surface $z=h(x,y, t)$ at each moment of time divides the region $\bar{\Omega}$ into two parts: $\bar{\Omega}^{(-)} =\{z:z\in[-a;h(x,y, t) ] \}$ and $\bar{\Omega}^{(+)}=\{z:z\in[h(x,y, t);a]\}$. Let us introduce the local coordinates $r,l,m$ and write the equation of the normal to the surface $h(l,m,t)$:
$$ \frac{x-l}{h_l (l,m,t)}=\frac{y-m}{h_m (l,m,t)}=\frac{z-h(l,m,t)}{-1}. $$
The equation for the distance from the point with coordinates $(r,l,m )$ to the surface $h(l,m,t)$ along the normal to it has the following form:
$$ r^2 = (x-l)^2 + (y-m)^2 + (z-h)^2 = (z-h)^2 (1+h_{l}^2+h_{m}^2).$$

For a detailed description of the solution in the inner transition layer, in the vicinity of surface $h(l,m,t)$ we proceed to the extended variable
\begin{equation} \label{xidefin}
\displaystyle \xi=\frac{r}{\mu}
\end{equation}
and to the local coordinates $(r,l,m)$ using the relations
\begin{equation} \label{localrelations}
x=l-\frac{r h_l}{ \sqrt{1+h_{l}^2+h_{m}^2}}, \quad y=m-\frac{r h_m}{ \sqrt{1+h_{l}^2+h_{m}^2}}, \quad z=h(l,m,t)+ \frac{r }{ \sqrt{1+h_{l}^2+h_{m}^2}}.
\end{equation}

We assume $r>0$ in the domain $ \mathbb{R}^2 \times \bar{\Omega}^{(+)}$ and $r<0$ in the domain $ \mathbb{R}^2 \times \bar{\Omega}^{(-)} $, and we note that if $z=h(l,m,t)$, then we have $r=0$, $l=x$, and $m=y$; the derivatives of the function $h(l,m,t)$ in \eqref{localrelations} are also taken for $l=x$ and $m=y$.

The asymptotic approximation $U(x,y,z,t,\mu)$ of the solution of problem \eqref{mainproblem} is constructed separately in regions $\bar{\Omega}^{(-)}$ and $\bar{\Omega}^{(+)}$, i.e.,
\begin{equation} \label{asymptoticsolution}
U(x, y, z,t, \mu)=
\begin{cases}
U^{(-)} (x, y, z,t, \mu) , \quad (x, y,z, t) \in \mathbb{R}^2\times \bar{\Omega}^{(-)} \times \mathcal{T} ,\\
U^{(+)}(x, y,z, t, \mu) , \quad (x, y,z,t )\in \mathbb{R}^2\times\bar{\Omega}^{(+)} \times \mathcal{T},
\end{cases}
\end{equation}
as the sums of two terms
\begin{equation} \label{asymptoticapproximation2terms}
U^{(\mp)}=\bar{u}^{(\mp)}(x, y,z, \mu)+Q^{(\mp)}(\xi, l,m, h(l,m, t), t, \mu),
\end{equation}
where $\bar{u}^{(\mp)}(x,y,z,t,\mu)$ are the functions describing the outer region and $Q^{(\mp)}(\xi,l,m,h(l,m,t),$ $t,\mu)$ are the functions describing the inner transition layer.

Each term in \eqref{asymptoticapproximation2terms} is represented as an expansion in powers of the small parameter $\mu$:
\begin{align}
\bar{u}^{(\mp)}(x, y,z, \mu)&=\bar{u}_{0}^{(\mp)}(x, y,z)+\mu \bar{u}_{1}^{(\mp)}(x, y,z)+\ldots, \label{expansionregularfunctions} \\
Q^{(\mp)}(\xi, l,m, h, t, \mu)&=Q_{0}^{(\mp)}(\xi, l,m, h, t)+\mu Q_{1}^{(\mp)}(\xi, l,m, h, t)+\ldots. \label{expansiontransitionfunctions}
\end{align}

We assume that $z=h(x,y,t)$ is the surface on which the solution $u(x,y,z,t,\mu)$ to problem \eqref{mainproblem} at each time $t$ takes on a value equal to the half-sum of the functions $\bar{u}^{(-)} (x, y,z)$ and $\bar{u}^{(+)}(x, y,z)$:
\begin{equation} \label{halfsum}
\displaystyle \phi(x,y, h(x,y, t), \mu):=\frac{1}{2}\left(\bar{u}^{(-)}(x,y, h(x, y,t),\mu)+\bar{u}^{(+)}(x,y, h(x, y,t),\mu) \right).
\end{equation}
For the zeroth-order approximation, it takes the form
\begin{equation}
\label{halfsum0order}
\displaystyle \phi_0(x,y, h(x,y, t)):=\frac{1}{2}\left(\varphi^{(-)}(x,y, h(x, y,t))+\varphi^{(+)}(x,y, h(x, y,t)) \right).
\end{equation}

The functions $U^{(-)}(x,y,z,t,\mu)$ and $U^{(+)}(x,y,z,t,\mu)$ and their derivatives along the normal to the surface $z=h(x,y,t)$ are matched continuously on the surface $h(x,y,t)$ at each moment of time $t$, i.e., the following two conditions hold:
\begin{align}
U^{(-)} (x,y, h(x,y, t), t, \mu)=U^{(+)}(x, y,h(x, y,t), t, \mu)=\phi(x,y, h(x,y, t), \mu) , \label{sewingcond1} \\
\displaystyle \frac{\partial U^{(-)}}{\partial n}(x,y, h(x,y, t), t, \mu)=\frac{\partial U^{(+)}}{\partial n}(x, y,h(x,y, t), t, \mu), \label{sewingcond2}
\end{align}
where the function $\phi(x,y, h(x,y, t), \mu)$ is defined in \eqref{halfsum}.

The surface $z=h(x,y,t)$ is also sought in the form of an expansion in powers of a small parameter, i.e.,
\begin{equation} \label{curveexpansion}
h(x,y, t)=h_{0}(x,y, t)+\mu h_{1}(x, y,t)+\mu^{2}h_{2}(x, y,t)+\ldots ,
\end{equation}
and we introduce a notation for the approximation of $h(x,y,t)$ with expansion terms up to order $n$, i.e.,
\begin{equation} \label{hExpansion}
\hat{h}_{n}(x,y,t)=\sum_{i=0}^{n} \mu^{i} h_{i}(x,y,t), \quad (x,y) \in \mathbb{R}^2, \ t \in \bar{\mathcal{T}}.
\end{equation}

We introduce the vectors ${\bf x}=(x,y,z,t)^{T}$ and ${\bf r}=(r,l,m,t)^{T}$, where $r=r(x,y,z,t), l=l(x,y,z,t), m=m(x,y,z,t)$, and $x,y,z$ are defined in \eqref{localrelations}.

The total differential of the vector ${\bf r}$ is $\mathrm{d}{\bf r}=\nabla_{{\bf x}}{\bf r} \ \mathrm{d}{\bf x}$. On the other hand, if we look at the inverse problem of determining the differential changes in our original coordinate system $(x,y,z,t)$ from those in the coordinate system $(r,l,m,t)$, we obtain $\mathrm{d}{\bf x}=\nabla_{{\bf r}}{\bf x} \ \mathrm{d} {\bf r}$, from which follows $\nabla_{{\bf x}} {\bf r}=\left(\nabla_{{\bf r}}{\bf x} \right)^{-1}$, which in matrix form is expressed as
\begin{align} \label{ComponentsOfNabla}
\left(\begin{smallmatrix} 
\displaystyle \frac{\partial r }{\partial x} & \displaystyle  \frac{\partial r }{\partial y} &\displaystyle \frac{\partial r }{\partial z} &\displaystyle \frac{\partial r }{\partial t}\\
\displaystyle \frac{\partial l }{\partial x} &\displaystyle \frac{\partial l }{\partial y} &\displaystyle \frac{\partial l }{\partial z}  &\displaystyle \frac{\partial l }{\partial t}\\
\displaystyle \frac{\partial m }{\partial x} &\displaystyle \frac{\partial m }{\partial y} &\displaystyle \frac{\partial m }{\partial z}&\displaystyle \frac{\partial m }{\partial t} \\
\displaystyle 0 &\displaystyle 0 &\displaystyle 0&\displaystyle 1
\end{smallmatrix} \right)= \left(\begin{smallmatrix} 
\displaystyle\frac{\partial x }{\partial r} &\displaystyle \frac{\partial x }{\partial l} &\displaystyle \frac{\partial x }{\partial m} &\displaystyle \frac{\partial x }{\partial t}\\\\
\displaystyle\frac{\partial y }{\partial r} &\displaystyle \frac{\partial y }{\partial l} &\displaystyle \frac{\partial y }{\partial m} &\displaystyle \frac{\partial y }{\partial t}\\
\displaystyle\frac{\partial z }{\partial r} &\displaystyle \frac{\partial z }{\partial l} &\displaystyle \frac{\partial z }{\partial m} &\displaystyle \frac{\partial z }{\partial t} \\
\displaystyle 0 &\displaystyle 0 &\displaystyle 0&\displaystyle 1
\end{smallmatrix} \right)^{-1}.
\end{align} 

We rewrite the differential operators in \eqref{mainproblem} in terms of the variables $r,l,m,t$:
\begin{multline} \label{nabla}
\displaystyle \nabla=\left\lbrace \frac{\partial }{\partial x}; \frac{\partial }{\partial y}; \frac{\partial }{\partial z} \right\rbrace=\left\lbrace \frac{\partial r }{\partial x} \frac{\partial }{\partial r}+\frac{\partial l }{\partial x} \frac{\partial }{\partial l}+\frac{\partial m }{\partial x} \frac{\partial }{\partial m} ;\frac{\partial r }{\partial y} \frac{\partial }{\partial r} +\frac{\partial l }{\partial y} \frac{\partial }{\partial l}+\frac{\partial m }{\partial y} \frac{\partial }{\partial m}; \right. \\
\left. \frac{\partial r }{\partial z} \frac{\partial }{\partial r} +\frac{\partial l }{\partial z} \frac{\partial }{\partial l}+\frac{\partial m }{\partial z} \frac{\partial }{\partial m} \right\rbrace ,
\end{multline}
where the components of the operator $\nabla$ can be found from \eqref{localrelations} and \eqref{ComponentsOfNabla}.

The differential operator $\Delta$ takes the form:
\begin{multline} \label{laplasian}
\Delta= \frac{\partial^2 }{\partial x^2}+ \frac{\partial^2 }{\partial y^2}+ \frac{\partial^2 }{\partial z^2} = \frac{\partial r }{\partial x} \frac{\partial }{\partial r}\left(\frac{\partial }{\partial x} \right)+\frac{\partial l }{\partial x} \frac{\partial }{\partial l}\left(\frac{\partial }{\partial x} \right)+\frac{\partial m }{\partial x} \frac{\partial }{\partial m}\left(\frac{\partial }{\partial x} \right) \\
+\frac{\partial r }{\partial y} \frac{\partial }{\partial r}\left(\frac{\partial }{\partial y} \right)+\frac{\partial l }{\partial y} \frac{\partial }{\partial l}\left(\frac{\partial }{\partial y} \right)+\frac{\partial m }{\partial y} \frac{\partial }{\partial m}\left(\frac{\partial }{\partial y} \right)\\
+\frac{\partial r }{\partial z} \frac{\partial }{\partial r}\left(\frac{\partial }{\partial z} \right)+\frac{\partial l }{\partial z} \frac{\partial }{\partial l}\left(\frac{\partial }{\partial z} \right)+\frac{\partial m }{\partial z} \frac{\partial }{\partial m}\left(\frac{\partial }{\partial z} \right),
\end{multline}
where the values for $ {\partial }/{\partial x}, {\partial }/{\partial y}, {\partial }/{\partial z} $ are taken from \eqref{nabla}.

We also rewrite the operator $\partial / \partial t$ in the variables $r,l,m,t$:
\begin{align} \label{operatordt}
\frac{\partial}{\partial t}= \frac{\partial}{\partial t} + \frac{\partial r}{\partial t}\frac{\partial }{\partial r}+\frac{\partial l}{\partial t} \frac{\partial }{\partial l} +\frac{\partial m}{\partial t} \frac{\partial }{\partial m},
\end{align} 
where components $\partial r/ \partial t, \partial l / \partial t,\partial m /\partial t $ are found from \eqref{ComponentsOfNabla}.

The derivative along the normal to the surface $h(x,y,t)$ with respect to the variables $x,y,z$ takes the form:
\begin{equation} \label{derivativetonormalXYZ}
\frac{\partial}{\partial n}=(\textbf{n}, \nabla)= -\frac{ h_x}{ \sqrt{1+h_{x}^2+h_{y}^2}}\frac{\partial}{\partial x}-\frac{ h_y}{ \sqrt{1+h_{x}^2+h_{y}^2}}\frac{\partial}{\partial y}+ \frac{1 }{ \sqrt{1+h_{x}^2+h_{y}^2}}\frac{\partial}{\partial z},
\end{equation}
where values for the vector $\textbf{n}$ are obtained from \eqref{localrelations}. With respect to the variables $\xi, l,m, t$, the derivative \eqref{derivativetonormalXYZ} has the form
\begin{equation} \label{derivativetonormal}
\frac{\partial}{\partial n}=\frac{\partial}{\partial r}=\frac{1}{\mu} \frac{\partial}{\partial \xi}.
\end{equation}

\subsection{Outer functions}

Let us find the outer functions of zeroth order by substituting the expansions \eqref{expansionregularfunctions} into the stationary equation:
\begin{equation}
\displaystyle \mu \Delta \bar{u} =-\bar{u} \left( \frac{\partial \bar{u}}{\partial x}+ \frac{\partial \bar{u}}{\partial y} + \frac{\partial \bar{u}}{\partial z}\right)+f(x, y,z).
\end{equation}
Expanding the functions into series in powers of a small parameter and equating the coefficients at the same power $\mu$, we obtain first-order PDEs for the functions $\bar{u}_{i}^{(\mp)}(x,y,z)$, $i=0,1$. More precisely, by equating the coefficients at $\mu^{0}$, we obtain
\begin{align}\label{zeroorderregularequation1}
\begin{cases}
\displaystyle \bar{u}_{0}^{(-)} \left(\frac{\partial\bar{u}_{0}^{(-)}}{\partial x}+\frac{\partial\bar{u}_{0}^{(-)}}{\partial y}+\frac{\partial\bar{u}_{0}^{(-)}}{\partial z}\right)=f(x, y,z),\\
\bar{u}_{0}^{(-)}(x, y,-a)=u^{-a}(x,y), \quad \bar{u}_{0}^{(-)}(x, y,z)=\bar{u}_{0}^{(-)}(x+L, y,z)=\bar{u}_{0}^{(-)}(x, y+M,z);
\end{cases}
\end{align}
\begin{align}\label{zeroorderregularequation2}
\begin{cases}
\displaystyle \bar{u}_{0}^{(+)}\left( \frac{\partial\bar{u}_{0}^{(+)}}{\partial x}+\frac{\partial\bar{u}_{0}^{(+)}}{\partial y}+\frac{\partial\bar{u}_{0}^{(+)}}{\partial z}\right)=f( x, y,z),\\
\bar{u}_{0}^{(+)}(x,y, a)=u^{a}(x,y), \quad \bar{u}_{0}^{(+)}(x, y,z)=\bar{u}_{0}^{(+)}(x+L, y,z)=\bar{u}_{0}^{(+)}(x, y+M,z).
\end{cases}
\end{align}

According to Assumption~\ref{A1}, \eqref{zeroorderregularequation1} and \eqref{zeroorderregularequation2} have solutions
\begin{align}
\label{u0regu}
\bar{u}_{0}(x,y,z)= \begin{cases}
\bar{u}_{0}^{(-)}(x,y,z)=\varphi^{(-)}(x,y,z), \quad (x,y,z)\in \mathbb{R}^2 \times \bar{\Omega}^{(-)},\\
\bar{u}_{0}^{(+)}(x,y,z)=\varphi^{(+)}(x,y,z), \quad (x,y,z)\in \mathbb{R}^2 \times \bar{\Omega}^{(+)},
\end{cases}
\end{align}
where $\varphi^{(\mp)}$ can be calculated explicitly; see \eqref{eq7}.

First-order approximations of outer functions $\bar{u}_{1}^{(\mp)}(x,y)$ are defined as solutions to the problems
\begin{equation} \label{firstorderregularequation}
\begin{dcases}
\frac{\partial\bar{u}_{1}^{(\mp)}}{\partial x} + \frac{\partial\bar{u}_{1}^{(\mp)}}{\partial y}+ \frac{\partial\bar{u}_{1}^{(\mp)}}{\partial z}+\bar{u}_{1}^{(\mp)} P^{(\mp)}(x,y,z)=W^{(\mp)}(x,y,z),\\
\bar{u}_{1}^{(-)}(x,y, -a)=0, \quad \bar{u}_{1}^{(+)}(x,y, a)=0,\\
\bar{u}_{0}^{(\mp)}(x, y,z)=\bar{u}_{0}^{(\mp)}(x+L, y,z)=\bar{u}_{0}^{(\mp)}(x, y+M,z),
\end{dcases}
\end{equation}
where
$$P^{(\mp)}(x,y,z)= \frac{1}{\varphi^{(\mp)}} \left( \frac{\partial\varphi^{(\mp)}}{\partial x} + \frac{\partial\varphi^{(\mp)}}{\partial y}+ \frac{\partial\varphi^{(\mp)}}{\partial z}\right), \quad W^{(\mp)}(x,y,z)= -\frac{\Delta \varphi^{(\mp)}}{\varphi^{(\mp)}} . $$

From \eqref{firstorderregularequation}, we obtain $\bar{u}_{1}^{(\mp)}$ in its explicit form:
\begin{align} \label{firstorderregularfunctionsExplicit}
\begin{dcases}
\bar{u}_{1}^{(-)}= \frac{\int _{ 0}^{z+a} \exp \left(\int _0^{\zeta+a }P^{(-)}\left({\scriptstyle x-z-a+s, y-z-a+s,-a+s }\right)ds\right) W^{(-)}\left({\scriptstyle x-z-a+\zeta , y-z-a+\zeta,-a+\zeta } \right)d\zeta }{\exp \left(\int _0^{z+a} P^{(-)}\left( {\scriptstyle x-z-a+s, y-z-a+s,-a+s }\right)ds\right)} , \\
\bar{u}_{1}^{(+)}= \frac{-\int _{z-a}^{0} \exp \left(\int _0^{\zeta-a }P^{(+)}\left({\scriptstyle x-z+a+s, y-z+a+s,a+s }\right)ds\right) W^{(+)}\left({\scriptstyle x-z+a+\zeta , y-z+a+\zeta,a+\zeta } \right)d\zeta }{\exp \left(\int _0^{z-a} P^{(+)}\left( {\scriptstyle x-z+a+s, y-z+a+s,a+s } \right)ds\right)} .
\end{dcases}
\end{align}

\subsection{Inner-layer functions}

Substituting \eqref{asymptoticapproximation2terms} into \eqref{mainproblem}, subtracting the outer part from the equation, and then moving to the variables $\xi,l,m,t$ using \eqref{nabla}, \eqref{laplasian}, and \eqref{operatordt}, we obtain the equations for the inner-transition-layer functions $Q^{(\mp)}(\xi, l,m, h(l,m, t), t, \mu)$:
\begin{multline} \label{transitionlayerequation}
\displaystyle \frac{1}{\mu} \Big(\frac{\partial^{2}Q^{(\mp)}}{\partial\xi^{2}}+\frac{h_{t}-(h_{l}+h_{m}-1)\left(\bar{u}^{(\mp)}(l, m,h(l,m,t), \mu) + Q^{(\mp)} \right) }{\sqrt{1+h_{l}^{2}+h_{m}^{2}}}\frac{\partial Q^{(\mp)}}{\partial\xi}\Big)\\ +\displaystyle \sum_{i=0}\mu^{i} L_{i}[Q^{(\mp)}]=0,
\end{multline}
where $L_{i}$ represents differential operators of first or second order in the variables $\xi,l,m,t$.

\subsubsection{Inner-layer functions of zeroth order}

The equation that determines the zeroth-order functions of the inner transition layer $Q_{0}^{(\mp)}(\xi,l,m,h_0,t)$ is found by substituting the series \eqref{expansionregularfunctions}, \eqref{expansiontransitionfunctions}, and \eqref{curveexpansion} into equations \eqref{sewingcond1} and \eqref{transitionlayerequation}, expanding all the terms of \eqref{transitionlayerequation} in series in powers of $\mu$, equating the coefficients at $\mu^{-1}$ in \eqref{transitionlayerequation} and at $\mu^{0}$ in \eqref{sewingcond1}, and taking into account an additional condition for the decay of the transition functions at infinity:
\begin{align} \label{transitionalfunczeroord}
\begin{dcases}
\displaystyle \frac{\partial^{2}Q_{0}^{(\mp)}}{\partial\xi^{2}}
+\displaystyle \frac{{h_0}_{t}+\left(1- {h_0}_{l}- {h_0}_{m}\right)\left( \varphi^{(\mp)}(l,m, h_0)+Q_{0}^{(\mp)} \right)}{\sqrt{1+{h_0}_{l}^{2}+{h_0}_{m}^{2}}}\frac{\partial Q_{0}^{(\mp)}}{\partial\xi}=0;\\
\varphi^{(\mp)}(l,m, h_0)+Q_{0}^{(\mp)}(0, l,m, h_0, t)=\phi_0(l,m, h_0),\quad Q_{0}^{(\mp)}(\mp\infty, l, m, h_0, t)=0.
\end{dcases}
\end{align}
Note that the function $Q_{0}^{(-)}$ is defined for $\xi\leq 0$, while the function $Q_{0}^{(+)}$ is defined for $\xi\geq 0$.

To find the function $h_0(l,m, t)$, we introduce the auxiliary function
\begin{equation} \label{auxilaryfunction}
\tilde{u}(\xi, h_0(l,m, t))=
\begin{cases}
\varphi^{(-)} (l,m, h_0)+Q_{0}^{(-)} (\xi, l,m, h_0, t) , \quad \xi\leq 0,\\
\varphi^{(+)}(l,m, h_0)+Q_{0}^{(+)}(\xi, l,m, h_0, t) , \quad \xi\geq 0,
\end{cases}
\end{equation}
and in this notation, each of the equations in \eqref{transitionalfunczeroord} takes the form
\begin{equation} \label{replacement}
\displaystyle \frac{\partial^{2}\tilde{u}}{\partial\xi^{2}}+\displaystyle \frac{{h_0}_{t}+\tilde{u} \left(1- {h_0}_{l}- {h_0}_{m}\right)}{\sqrt{1+{h_0}_{l}^{2}+{h_0}_{m}^{2}}}\frac{\partial\tilde{u}}{\partial\xi}=0,
\end{equation}
where the variables $l,m,t$ and the function $h(l,m,t)$ are parameters.

Denoting $\displaystyle \frac{\partial \tilde{u} }{\partial\xi} = g(\tilde{u})$ and $\frac{\partial^2 \tilde{u} }{\partial\xi^{2}} = \frac{\partial g(\tilde{u}) }{\partial \tilde{u} } g(\tilde{u})$, \eqref{replacement} is transformed to
$$\displaystyle \frac{\partial g(\tilde{u}) }{\partial \tilde{u} } = \displaystyle \frac{-{h_0}_{t}+\tilde{u} \left({h_0}_{l}+ {h_0}_{m}-1 \right)}{\sqrt{1+{h_0}_{l}^{2}+{h_0}_{m}^{2}}} ,$$
from which we can deduce that
\begin{align} \label{derivativetildeu}
\frac{\partial\tilde{u}}{\partial\xi}= \begin{cases}
\displaystyle \Phi^{(-)} (\xi, h_0(l,m, t))=\int_{\varphi^{(-)}}^{\tilde{u}} \displaystyle \frac{-{h_0}_{t}+\tilde{u} \left({h_0}_{l}+ {h_0}_{m}-1 \right)}{\sqrt{1+{h_0}_{l}^{2}+{h_0}_{m}^{2}}} du , \quad \xi\leq 0,\\
\displaystyle \Phi^{(+)}(\xi, h_0(l,m, t))= \int_{\varphi^{(+)}}^{\tilde{u}} \displaystyle \frac{-{h_0}_{t}+\tilde{u} \left({h_0}_{l}+ {h_0}_{m}-1 \right)}{\sqrt{1+{h_0}_{l}^{2}+{h_0}_{m}^{2}}} du , \quad \xi\geq 0.
\end{cases}
\end{align}

Considering expansions \eqref{asymptoticapproximation2terms}--\eqref{curveexpansion} and \eqref{derivativetonormal} and equating the coefficients at $\mu^{-1}$ in \eqref{sewingcond2}, we obtain [given that at $\xi=0$, we have $x=l$, $y=m$, and $z=h_0(x,y,t)$]
\begin{align}\label{sewindcondexpanded0}
\Phi^{(-)} (0, h_0(x,y, t))-\Phi^{(+)} (0, h_0(x,y, t))=0.
\end{align}
From \eqref{derivativetildeu} and \eqref{sewindcondexpanded0}, we obtain
\begin{equation*}
\int_{\varphi^{(-)} (x, y,h_0(x,y, t))}^{\varphi^{(+)}(x,y, h_0(x,y, t))} \displaystyle \frac{-{h_0}_{t}+\tilde{u} \left({h_0}_{x}+ {h_0}_{y}-1 \right)}{\sqrt{1+{h_0}_{x}^{2}+{h_0}_{y}^{2}}} du =0,
\end{equation*}
and the equation that determines the zeroth-order approximation of the surface $h_0(x,y,t)$ takes the form
\begin{equation} \label{h0hxhyequation}
{h_0}_{t}=\frac{1}{2} \left({h_0}_{x}+{h_0}_{y}-1\right) \left(\varphi^{(+)}(x,y, h_0)+\varphi^{(-)}(x,y, h_0) \right)
\end{equation}
with the additional boundary condition ${h_0}(x,y,0)=h_{0}^{*} \in \bar{\Omega}$ and periodic conditions on the $x$ and $y$ axes. By virtue of Assumptions~\ref{A2} and \ref{A3}, equation \eqref{h0hxhyequation} has a solution inside the region $\Omega$ for any $(x,y,t) \in \mathbb{R}^2\times \bar{\mathcal{T}}$.

We solve \eqref{derivativetildeu}, and with \eqref{h0hxhyequation}, the functions $Q_{0}^{(\mp)}(\xi, l,m, h_0, t)$ are given by expressions in which $h_0(x,y,t)$ is a parameter:
\begin{align} \label{Q0equation}
\displaystyle Q_{0}^{(\mp)}(\displaystyle \xi,l,m, h_0, t)=\frac{2P^{(\mp)}(x,y,h_0)}{\exp \left(-\xi \frac{ P^{(\mp)}(x,y,h_0) \left(1-{h_0}_{x}-{h_0}_{y} \right)}{\sqrt{1+{h_0}_{x}^{2}+{h_0}_{y}^{2}}} \right)+1},
\end{align}
where
\begin{align*}
P^{(-)}(x,y,h_0)=\frac{1}{2}\left(\varphi^{(+)}(x,y,h_0)-\varphi^{(-)}(x,y,h_0) \right),\\
P^{(+)}(x,y,h_0)=\frac{1}{2}\left(\varphi^{(-)}(x,y,h_0)-\varphi^{(+)}(x,y,h_0)\right).
\end{align*}

The transition-layer functions $Q_{0}^{(\mp)}(\displaystyle \xi,l,m, h_0, t)$ are exponentially decreasing with $\xi \rightarrow \mp \infty$ and satisfy the exponential estimates \cite{Vasileva1998ContrastSI,Butuzov1997ASYMPTOTICTO}
\begin{equation}\label{equat22}
\underaccent{\bar}{C} e^{\underaccent{\bar}{\kappa} \xi} \leq |Q_{0}^{(-)}( \xi, l,m, h_0, t)|\leq \bar{C}e^{\bar{\kappa}\xi} , \quad \xi\leq 0, \quad t\in \bar{\mathcal{T}},
\end{equation}
\begin{equation}\label{equat23}
\underaccent{\bar}{C} e^{-\underaccent{\bar}{\kappa} \xi} \leq |Q_{0}^{(+)}( \xi, l,m, h_0, t)|\leq \bar{C}e^{-\bar{\kappa}\xi}, \quad \xi\geq 0, \quad t\in \bar{\mathcal{T}},
\end{equation}
where $\underaccent{\bar}{C}, \bar{C}$ and $\underaccent{\bar}{\kappa} , \bar{\kappa}$ are four positive constants that are independent of $\xi,l,m,t$; in particular, we have $$\underaccent{\bar}{C} := \frac{1}{2} \inf_{t \in \bar{\mathcal{T}}} \varphi^{(+)}(x,y,h_{0})- \varphi^{(-)}(x,y,h_{0}), $$
$$ \bar{C} := \frac{1}{2} \sup_{t \in \bar{\mathcal{T}}} \varphi^{(+)}(x,y,h_{0})- \varphi^{(-)}(x,y,h_{0}).$$

From the boundary conditions of \eqref{transitionalfunczeroord}, we deduce that
$$|Q_{0}^{( \mp)}( 0, x,y, h_0, t)|=|\phi_0(x, y,h_0)-\varphi^{(\mp)}(x,y, h_0)| > \frac{1}{2} \left(u^{a}(x,y)-u^{-a}(x,y) \right) > \mu^2 $$
and $|Q_{0}^{(\mp)}( \xi, l,m, h_0, t)|\rightarrow 0$ for $\xi \rightarrow \mp \infty$. Because $\mu >0$ is a constant, $|Q_{0}^{(\mp)}( \xi, l,m, h_0, t)|$ is a decreasing function and $\xi= \frac{1}{\mu} \sqrt{1+h_{l}^2+h_{m}^2} \left(z- h(l,m, t)\right) $, then there exist $H^{(\mp)}(l,m, t)$ such that on the intervals $[-a,H^{(-)}(l,m, t)]$ and $[H^{(+)}(l,m, t),a]$ we have $|Q_{0}^{(\mp)}(\xi, l, $ $m, h_0, t)| \leq \mu^2 $ for every $l,m,t$, and at $z=H^{(\mp)}$, we have
\begin{equation}
\label{Qlr-mu}
|Q_{0}^{(-)}(\xi( H^{(-)}), l,m, H^{(-)}, t)| = \mu^2, \quad |Q_{0}^{(+)}(\xi( H^{(+)}), l,m, H^{(+)}, t)| = \mu^2.
\end{equation}

We denote the width of the transition layer as
\begin{equation*}
\Delta h = \max \limits_{\begin{subarray}{c} (x,y) \in \mathbb{R}^2 \\ t\in \bar{\mathcal{T}} \end{subarray} } {\left( H^{(+)}(x,y, t) -H^{(-)}(x,y, t) \right)},
\end{equation*}
and from \eqref{equat22} and \eqref{equat23} for $z=H^{(\mp)}$, it follows that
\begin{equation}
\displaystyle \underaccent{\bar}{C} e^{- \frac{\underaccent{\bar}{\kappa}}{2\mu } \sqrt{1+h_{l}^2+h_{m}^2} \Delta h} \leq \mu^2 \leq \bar{C}e^{- \frac{\bar{\kappa}}{2\mu } \sqrt{1+h_{l}^2+h_{m}^2} \Delta h}.
\end{equation}
Thus, the width of the inner transition layer $\Delta h$ can be estimated as
\begin{equation} \label{Deltah}
\frac{2\mu }{\underaccent{\bar}{\kappa} \sqrt{1+h_{l}^2+h_{m}^2}} \ln \frac{\underaccent{\bar}{C}}{\mu^2} \leq \Delta h \leq \frac{2\mu }{\bar{\kappa} \sqrt{1+h_{l}^2+h_{m}^2}} \ln \frac{\bar{C}}{\mu^2}, \text{~i.e.~} \Delta h \sim \mu |\ln \mu|.
\end{equation}

\subsubsection{Inner-layer functions of first order}

The equation for the first-order transition-layer functions $Q_{1}^{(\mp)}(\xi,l,m,\hat{h}_1,t)$ is obtained by equating the terms containing $\mu^{0}$ in \eqref{transitionlayerequation}:
\begin{multline}\label{transitionlayerfirstord}
\displaystyle \frac{\partial^{2}Q_{1}^{(\mp)}}{\partial\xi^{2}}+ \frac{\partial}{\partial \xi} \left(Q_{1}^{(\mp)} \frac{{h_0}_{t}  + \tilde{u}(1-{h_0}_{l}-{h_0}_{m} ) }{ \sqrt{1+{h_0}_{l}^{2}+{h_0}_{m}^{2} }} \right) \\ = r_{1}^{(\mp)}(\xi,l,m,t) {h_1}_{l}
+r_{2}^{(\mp)}(\xi,l,m,t) {h_1}_{m} +r_{3}^{(\mp)}(\xi,l,m,t) h_{1} +r_{4}^{(\mp)}(\xi,l,m,t){h_1}_{t}+r_{5}^{(\mp)}(\xi,l,m,t) \\ :=f_{1}^{(\mp)}(\xi,l,m,t),
\end{multline}
where $r_{1}^{(\mp)}(\xi,l,m,t)$ to $r_{5}^{(\mp)}(\xi,l,m,t)$ are known functions, and in particular, $r_{4}^{(\mp)}(\xi,l,m,t)=\frac{\partial Q_{0}^{(\mp)} (\xi,l,m, h_0,t)}{\partial \xi} \frac{1}{\sqrt{1+{h_0}_{l}^{2}+{h_0}_{m}^{2}}}$. The derivatives of the function $ h_1 (l,m, t) $ are taken for $ l=x $ and $ m=y $. Equation~\eqref{sewingcond1} at order $\mu^{1}$ implies the boundary conditions
\begin{align}\label{transitionlayerfirstordbound1}
Q_{1}^{(\mp)}(0, l,m, \hat{h}_1, t)=\frac{1}{2}\left( \bar{u}_{1}^{(\pm)}(l,m, h_0)-\bar{u}_{1}^{(\mp)}(l,m, h_0) \right) := p_1^{(\mp)}(l,m, h_0).
\end{align}
We also add conditions at infinity:
\begin{equation}\label{transitionlayerfirstordbound2}
Q_{1}^{(\mp)} (\mp \infty, l,m, \hat{h}_1, t)=0.
\end{equation}
Solutions to problems \eqref{transitionlayerfirstord}--\eqref{transitionlayerfirstordbound2} can be written explicitly as
\begin{align} \label{Q1function}
\begin{split}
Q_{1}^{(\mp)}(\xi, l,m, \hat{h}_1, t)=J^{(\mp)}(\xi,h_0) \left( p_{1}^{(\mp)}(l,m, h_0) \vphantom{\int_{0}^{0}} \right. \\
\left. +\int_{0}^{\xi}\frac{ds}{J^{(\mp)}(s,h_0)}\int_{\mp\infty}^{s}f_{1}^{(\mp)}(\eta, l,m, t)d\eta ds \right),
\end{split}
\end{align}
where $\displaystyle J^{(\mp)}(\xi,h_0)= \left( \Phi^{(\mp)}(0,h_0) \right)^{-1} \Phi^{(\mp)}(\xi,h_0)$.

Clearly, $Q_{1}^{(\mp)}(\xi, l,m, \hat{h}_1, t)$ satisfy exponential estimates \eqref{equat22} and \eqref{equat23}. From \eqref{Q1function}, we find
\begin{align} \label{Q1functionDerivative}
\begin{split}
\frac{\partial Q_{1}^{( \mp)}}{\partial \xi} (0, x,y, \hat{h}_1, t)=p_{1}^{(\mp)}(x,y, h_0) \left( -\frac{{h_0}_{t}}{\sqrt{1+{h_0}_{x}^{2}+{h_0}_{y}^{2}}} \vphantom{\frac{\phi(x,y, h_0)\left({h_0}_{x}+{h_0}_{y}-1\right)}{\sqrt{1+{h_0}_{x}^{2}+{h_0}_{y}^{2}}}} \right. \\
\left. +\frac{\phi_0(x,y, h_0)\left({h_0}_{x}+{h_0}_{y}-1\right)}{\sqrt{1+{h_0}_{x}^{2}+{h_0}_{y}^{2}}} \right) - \int_{0}^{\mp \infty} f_{1}^{( \mp)} (\eta,x,y,t) d\eta.
\end{split}
\end{align}

From the first-order $C^1$-matching condition \eqref{sewingcond2} and considering expansions \eqref{derivativetonormalXYZ} and \eqref{derivativetonormal}, we obtain
\begin{multline} \label{matchingfirstord}
h_1(x,y, t) \left(\frac{\partial^2 Q_{0}^{(-)}}{\partial \xi \partial h_0}-\frac{\partial^2 Q_{0}^{(+)}}{\partial \xi \partial h_0} \right)+\frac{\partial Q_{1}^{(-)}}{\partial \xi }-\frac{\partial Q_{1}^{(+)}}{\partial \xi }\\
+\frac{ {h_0}_{y} \left( \frac{\partial \varphi^{(+)}}{\partial y}-\frac{\partial \varphi^{(-)}}{\partial y} \right)
+{h_0}_{x} \left( \frac{\partial \varphi^{(+)}}{\partial x}-\frac{\partial \varphi^{(-)}}{\partial x} \right)
+\left( \frac{\partial \varphi^{(-)}}{\partial h_0}-\frac{\partial \varphi^{(+)}}{\partial h_0} \right)}{\sqrt{1+{h_0}_{x}^{2}+{h_0}_{y}^{2}}} =0.
\end{multline}
From \eqref{transitionlayerfirstord}, \eqref{Q1functionDerivative}, and \eqref{matchingfirstord}, we derive the equation that determines $h_1(x,y,t)$:
\begin{align}\label{eqforh1}
\begin{dcases}
{h_1}_{t} \frac{ \varphi^{(+)}-\varphi^{(-)} }{\sqrt{1+{h_0}_{x}^{2}+{h_0}_{y}^{2}}} ={h_1}_{x} V_1(x,y,t)+{h_1}_{y}V_2(x,y,t) + h_1 V_3(x,y,t) +V_4(x,y,t), \\
h_{1}(x,y, t)=h_{1}(x+L,y, t)=h_{1}(x,y+M, t),\ h_1(x,y,0)=0,
\end{dcases}
\end{align}
where $V_1(x,y,t)$ to $V_4(x,y,t)$ are known functions. The problem for the function $h_{1}(x,y, t)$ is solvable because the coefficient of the term ${h_1}_{t}$ in \eqref{eqforh1} is positive. Similarly to \eqref{firstorderregularequation}--\eqref{eqforh1}, it is possible to obtain asymptotic approximation terms up to order $n$.

\begin{longtable}{l l l}
\hline \multicolumn{1}{c}{\textbf{Notation}} & \multicolumn{1}{c}{\textbf{Description }} & \multicolumn{1}{c}{\textbf{Reference}} \\ \hline
\endfirsthead
\hline
\hspace{-3.5mm} \begin{tabular}{l} $\mu$ \end{tabular} & \hspace{-3.5mm} \begin{tabular}{l} Small parameter, $ 0<\mu \ll 1 $ \end{tabular} & Eq.\eqref{mainproblem} \\
\hline
\renewcommand{\arraystretch}{1.2} \hspace{-4.5mm} \begin{tabular}{l} $\Omega$, $\bar{\mathcal{T}}$ \end{tabular} & \hspace{-3.5mm} \begin{tabular}{l} Domains: $\Omega=(-a,a)$, $\bar{\mathcal{T}}=[0, T]$ \end{tabular} & Eq.\eqref{mainproblem} \\
\hline
\hspace{-3.5mm} \begin{tabular}{l} $u^{(\mp)}$ \end{tabular} & \hspace{-3.5mm} \begin{tabular}{l} Left and right boundary conditions \end{tabular} & Eq.\eqref{mainproblem} \\
\hline
\hspace{-3.5mm} \begin{tabular}{l} $l,m,r$ \end{tabular} & \hspace{-3.5mm} \begin{tabular}{l} The local coordinates which depend on $x,y,z$ \end{tabular} & Eq.\eqref{localrelations} \\
\hline
\hspace{-3.5mm} \begin{tabular}{l} $h(x,y,t)$ \end{tabular} & \hspace{-3.5mm} \begin{tabular}{l} The surface that is located in the middle \\ of the transition layer \end{tabular} & Eq.\eqref{curveexpansion} \\
\hline
\hspace{-3.5mm} \begin{tabular}{l} $\hat{h}_n(x,y,t)$ \end{tabular} & \hspace{-3.5mm} \begin{tabular}{l} The expansion of $h(x,y,t)$ in a series \\ up to order $n$ \end{tabular} & Eq.\eqref{hExpansion} \\
\hline
\hspace{-3.5mm} \begin{tabular}{l} $\displaystyle \frac{\partial }{\partial n}$ \end{tabular} & \hspace{-3.5mm} \begin{tabular}{l} Derivative along the normal \\ to the surface $h(x,y,t)$ \end{tabular} & Eq.\eqref{derivativetonormal} \\
\hline
\hspace{-3.5mm} \begin{tabular}{l} Superscript $ ^{(\mp)}$ \end{tabular} & \hspace{-3.5mm} \begin{tabular}{l} Describe functions on the left and right, \\
respectively, relative to the surface $h(x,y,t)$
\end{tabular} & Eq.\eqref{asymptoticsolution} \\
\hline
\hspace{-3.5mm} \begin{tabular}{l} Subscript $ _{0,1...}$ \end{tabular} & \hspace{-3.5mm} \begin{tabular}{l} The order of approximation of the \\ asymptotic solution \end{tabular} & Eq.\eqref{expansionregularfunctions}-\eqref{curveexpansion} \\
\hline
\hspace{-4.5mm} \renewcommand{\arraystretch}{1.2} \begin{tabular}{l} $\bar{\Omega}^{(-)}$, $\bar{\Omega}^{(+)}$ \end{tabular} & \hspace{-3.5mm} \begin{tabular}{l} $[-a,h(x,y,t)]$ and $[h(x,y,t),a]$ respectively \end{tabular} & Eq.\eqref{asymptoticsolution} \\
\hline
\hspace{-3.5mm} \begin{tabular}{l} $\bar{u}^{(\mp)}$ \end{tabular} & \hspace{-3.5mm} \begin{tabular}{l} Outer functions describing the \\ solution far from the surface $h(x,y,t)$ \end{tabular} & Eq.\eqref{asymptoticapproximation2terms} \\
\hline
\hspace{-3.5mm} \begin{tabular}{l} $\varphi^{(\mp)}$ \end{tabular} & \hspace{-3.5mm} \begin{tabular}{l} Zero-order approximation \\ of outer functions \end{tabular} & Eq.\eqref{eq7},\eqref{u0regu} \\
\hline
\hspace{-3.5mm} \begin{tabular}{l} $Q^{(\mp)}$ \end{tabular} & \hspace{-3.5mm} \begin{tabular}{l} Transition-layer functions describing \\ the solution near the surface $h(x,y,t)$ \end{tabular} & Eq.\eqref{expansiontransitionfunctions} \\
\hline
\hspace{-3.5mm} \begin{tabular}{l} $ \Delta h$ \end{tabular} & \hspace{-3.5mm} \begin{tabular}{l} The width of the transition layer, \\ $ \Delta h \sim \mu |\ln \mu|$ \end{tabular} & Eq.\eqref{Deltah} \\
\hline
\hspace{-3.5mm} \begin{tabular}{l} $ \xi$ \end{tabular} & \hspace{-3.5mm} \begin{tabular}{l} Extended variable, $\xi = r/ \mu$ \end{tabular} & Eq.\eqref{xidefin} \\
\hline
\hspace{-4.5mm} \renewcommand{\arraystretch}{1.2}
\begin{tabular}{l} $u^\varepsilon (x,y,z,t)$ \end{tabular} & \hspace{-3.5mm} \begin{tabular}{l} The smooth approximate data \end{tabular} & Eq.\eqref{uAlphaL},\eqref{uAlphaR} \\
\hline
\hspace{-4.5mm} \renewcommand{\arraystretch}{1.2}
\begin{tabular}{l} $f^*(x,y,z) $ \end{tabular} & \hspace{-3.5mm} \begin{tabular}{l} The exact source function \end{tabular} & Eq.\eqref{f0Ineq} \\
\hline
\hspace{-4.5mm} \renewcommand{\arraystretch}{1.2}
\begin{tabular}{l} $f^\delta(x,y,z) $ \end{tabular} & \hspace{-3.5mm} \begin{tabular}{l} The regularized approximate \\ source function \end{tabular} & Eq.\eqref{fdelta} \\
\hline
\caption{Notations and references to their definitions.}
\label{NotationTable}
\end{longtable}

\section{Main results} \label{statementresults}

For readability, in Table~\ref{NotationTable} we summarize the notations and abbreviations that are used frequently in this section. Below, we present the main results of this paper, which are based on the results of the asymptotic analysis from Section~\ref{constructionOfAsymptotic}.

\begin{theorem} \label{MainThm}
Suppose that functions $f(x,y,z), u^{-a}(x,y), u^{a}(x,y), u_{init}(x,y,z)$ are sufficiently smooth, $x$- and $y$-periodic, and $\mu\ll 1$. Then, under Assumptions~\ref{A1}--\ref{A4}, the boundary-value problem \eqref{mainproblem} has a unique smooth solution with an inner transition layer. In addition, the order-$n$ asymptotic solution $U_{n}(x,y,z,t,\mu)$ has the following representation:
\begin{align} \label{asymptoticnorder}
U_{n}= \begin{cases}
\displaystyle U_{n}^{(-)}=\sum_{i=0}^{n} \mu^{i} \left(\bar{u}_{i}^{(-)}\left(x,y,z\right)+Q_{i}^{(-)}\left(\xi_{i},l,m,\hat{h}_{i},t \right) \right), \ (x,y,z,t) \in \mathbb{R}^2 \times \bar{\Omega}^{(-)} \times \bar{\mathcal{T}} , \\
\displaystyle U_{n}^{(+)}=\sum_{i=0}^{n} \mu^{i} \left(\bar{u}_{i}^{(+)}\left(x,y,z\right)+Q_{i}^{(+)}\left(\xi_{i},l,m,\hat{h}_{i},t\right)\right), \ (x,y,z,t) \in \mathbb{R}^2 \times \bar{\Omega}^{(+)}\times \bar{\mathcal{T}},
\end{cases}
\end{align}
where $\xi_{i} = \left(z-\hat{h}_i \right) \sqrt{1+(\hat{h}_i)_{x}^{2}+(\hat{h}_i)_{y}^{2}}/ \mu$.

Moreover, the following asymptotic estimates hold:
\begin{equation} \label{NorderEstim1}
\forall (x,y,z,t)\in \mathbb{R}^2 \times \bar{\Omega} \times \bar{\mathcal{T}}:~ |u(x,y,z,t)-U_{n}(x,y,z,t,\mu)|=\mathcal{O}(\mu^{n+1}) ,
\end{equation}
\begin{equation} \label{NorderEstim2}
\forall (x,y,t) \in \mathbb{R}^2 \times \bar{\mathcal{T}}:~ |h(x,y,t)-\hat{h}_{n}(x,y,t)|=\mathcal{O}(\mu^{n+1}),
\end{equation}
\begin{equation}\label{NorderEstim3}
\forall (x,y,z,t) \in \mathbb{R}^2 \times \bar{ \Omega} \backslash \{ \hat{h}_{n}(x,y,t) \} \times \bar{\mathcal{T}}:~ \left| \frac{ \partial u(x,y,z,t)}{ \partial n}-\frac{ \partial U_{n}(x,y,z,t,\mu) }{ \partial n} \right| =\mathcal{O}(\mu^{n}).
\end{equation}
\end{theorem}

\begin{corollary}
\label{Corollary1}
(Zeroth-order approximation) Under the assumptions of Theorem~\ref{MainThm}, the zeroth-order asymptotic solution $U_{0}$ has the following representation:
\begin{align} \label{eq001}
U_{0}(x,y,z,t)=\begin{cases}
\varphi^{(-)}(x,y,z)+Q_{0}^{(-)}(\xi_{0},l,m,h_{0}(l,m,t),t) , \ (x,y,z,t)\in \mathbb{R}^2 \times\bar{\Omega}^{(-)} \times \bar{\mathcal{T}},\\
\varphi^{(+)}(x,y,z)+Q_{0}^{(+)}(\xi_{0},l,m,h_{0}(l,m,t),t) , \ (x,y,z,t)\in \mathbb{R}^2 \times\bar{\Omega}^{(+)} \times \bar{\mathcal{T}},
\end{cases}
\end{align}
where $\xi_{0}=\left(z-h_0 \right) \sqrt{1+(h_0)_{x}^{2}+(h_0)_{y}^{2}}/ \mu$. Moreover, the following hold:
\begin{equation} \label{eq002}
\forall (x,y,z,t)\in \mathbb{R}^2 \times \bar{\Omega} \times \bar{\mathcal{T}}:~ |u(x,y,z,t)-U_{0}(x,y,z,t)|=\mathcal{O}(\mu) ,
\end{equation}
\begin{equation} \label{eq003}
\forall t \in \bar{\mathcal{T}}:~ |h(x,y,t)-h_{0}(x,y,t)|=\mathcal{O}(\mu) .
\end{equation}

Furthermore, outside the narrow region $\left(h_{0}\left(x,y,t\right)- \Delta h /2, h_{0}(x,y,t)+\Delta h /2\right)$ with $\Delta h\sim\mu|\ln\mu|$, there exists a constant $C$ independent of $x, y, z, t,\mu$ such that the following inequalities hold:
\begin{equation} \label{eq004}
|u(x,y,z,t)-\varphi^{(-)}(x,y,z)| \leq C \mu, \qquad (x,y,z,t) \in \mathbb{R}^2 \times [-a, h_{0}(x,y,t)- \Delta h /2] \times \bar{\mathcal{T}},
\end{equation}
\begin{equation} \label{eq005}
|u(x,y,z,t)-\varphi^{(+)}(x,y,z)| \leq C \mu, \qquad (x,y,z,t) \in \mathbb{R}^2 \times [h_{0}(x,y,t)+\Delta h /2, a] \times \bar{\mathcal{T}},
\end{equation}
\begin{equation}\label{0orderEstim1}
\left| \frac{\partial u(x,y,z,t)}{\partial n}-\frac{\partial \varphi^{(-)}(x,y,z)}{\partial n} \right| \leq C \mu, \quad (x,y,z,t) \in \mathbb{R}^2 \times [-a, h_{0}(x,y,t)- \Delta h /2] \times \bar{\mathcal{T}},
\end{equation}
\begin{equation}\label{0orderEstim2}
\left| \frac{\partial u(x,y,z,t)}{\partial n}-\frac{\partial \varphi^{(+)}(x,y,z)}{\partial n} \right| \leq C \mu, \quad (x,y,z,t) \in \mathbb{R}^2 \times [h_{0}(x,y,t)+\Delta h /2, a] \times \bar{\mathcal{T}}.
\end{equation}
\end{corollary}

Corollary~\ref{Corollary1} follows directly from Theorem~\ref{MainThm}. Inequalities \eqref{0orderEstim1} and \eqref{0orderEstim2} in Corollary~\ref{Corollary1} can be obtained by taking into account the fact that the transition-layer functions are decreasing functions with respect to $\xi_0$ and are sufficiently small at the boundaries of the narrow region $(h_0(x,y,t)-\Delta h/2, h_0(x,y,t)+\Delta h/2)$, i.e., equation \eqref{Qlr-mu}.

From Corollary~\ref{Corollary1}, it follows that the solution can be approximated by outer functions of zeroth order everywhere except for a thin transition layer. In view of this, we construct the approximate source function $f$ by using only the outer functions of zeroth order, which reduces the computational cost significantly while retaining the accuracy of the results in the case of a relatively high noise level and small $\mu$. To do this, suppose that we have the deterministic noise model
\begin{multline}
\label{noisyData1}
\max \limits_{i,j,k } \left( \lvert u(x_i,y_j,z_k,t_0)-u^{\delta}_{i,j,k} \rvert, \lvert u_x(x_i,y_j,z_k,t_0)-(u_x)^{\delta}_{i,j,k} \rvert, \right. \\ \left. \lvert u_y(x_i,y_j,z_k,t_0)-(u_y)^{\delta}_{i,j,k} \rvert, \lvert u_z(x_i,y_j,z_k,t_0)-(u_z)^{\delta}_{i,j,k} \rvert \right) \leq \delta
\end{multline}
between the noisy data $\{u^{\delta}_{i,j,k}, (u_x)^{\delta}_{i,j,k}, (u_y)^{\delta}_{i,j,k}, (u_z)^{\delta}_{i,j,k} \}$ and the corresponding exact quantities $\{u(x_i,y_j,z_k,t_0), u_y(x_i,y_j,z_k,t_0),u_x(x_i,y_j,z_k,t_0),u_z(x_i,y_j,z_k,t_0) \}$ at time $t_0$ and at grid points $X_p:= \{x_0 < x_1 < \cdots < x_p \} $, $Y_q := \{ y_0 < y_1 < \cdots < y_q \}$, and $Z_v := \{-a=y_0 < y_1 < \cdots < y_v=a \}$, with maximum mesh sizes in each direction: $d_1:= \max\limits_{ i\in \{0, \cdots,p-1\} } \{x_{i+1} - x_{i}\}$, $d_2:= \max\limits_{ j\in \{0, \cdots,q-1\} } \{ y_{j+1} - y_{j}\}$, and $d_3:= \max\limits_{ k\in \{0, \cdots,v-1\} } \{ z_{k+1} - z_{k}\}$.

To study the convergence, we define the pre-approximated source function $f_0$ as
\begin{equation} \label{f0}
f_0(x,y,z) = u (x,y,z,t_0) \left( \frac{\partial u (x,y,z,t_0)}{\partial x} +\frac{\partial u (x,y,z,t_0)}{\partial y} +\frac{\partial u (x,y,z,t_0)}{\partial z} \right).
\end{equation}

\begin{proposition}
\label{ProAsympErr}
Let $f^*$ be the exact source function satisfying the original governing equation \eqref{mainproblem}. Under Assumption~\ref{A2}, there exists a constant $C_1$ independent of $x,y,\mu$ such that for any $p\in(0,+\infty)$,
\begin{equation}
\label{f0Ineq}
\|f^* - f_0 \|_{L^p(\mathbb{R}^2 \times\Omega)} \leq C_1 \mu |\ln \mu | .
\end{equation}
\end{proposition}

According to Corollary~\ref{Corollary1}, we can exclude data values $\{u^{\delta}_{i,j,k}, (u_x)^{\delta}_{i,j,k}, (u_y)^{\delta}_{i,j,k}, (u_z)^{\delta}_{i,j,k} \}$ belonging to the inner transition layer and use only nodes from the two regions $[-a, h_0 (x,y,t)$ $-\Delta h/2]$ and $[h_0 (x,y,t)+\Delta h/2,a]$ with node indices $k\in \mathcal{K}$, where $\mathcal{K} := \{ 0, \cdots, v^{(-)},v^{(+)},$ $ \cdots, v \}$. We also introduce notations for the left region, $\bar{\Omega}^l \equiv [-a, -a+d_3 v^{(-)}]$, and for the right region, $\bar{\Omega}^r \equiv [-a+d_3 v^{(+)},a]$.

We restore the source function $f^\delta(x,y,z)$ by using the following nonlinear formula:
\begin{equation}
\label{fdelta}
f^\delta = u^\delta (x,y,z,t_0) \left( \frac{\partial u^\delta (x,y,z,t_0)}{\partial x} +\frac{\partial u^\delta (x,y,z,t_0)}{\partial y} +\frac{\partial u^\delta (x,y,z,t_0)}{\partial z} \right).
\end{equation}

If we know only the measurements $u^\delta_{i,j,k}$, then we replace the values $\{u^{\delta}_{i,j,k}, ({u_x})^{\delta}_{i,j,k}, $ $ ({u_y})^{\delta}_{i,j,k}, ({u_z})^{\delta}_{i,j,k} \}$ in \eqref{fdelta} with the values $\{u^\varepsilon, \frac{\partial u^\varepsilon}{\partial x}, \frac{\partial u^\varepsilon}{\partial y}, \frac{\partial u^\varepsilon}{\partial z}\}$, where the function $u^\varepsilon (x,y,z,t_0)$ is constructed according to the following optimization problem for the left and right regions, respectively:
\begin{multline}
\label{uAlphaL}
u^\varepsilon = \arg\min_{\begin{subarray}{c} w\in C^1(\mathbb{R}^2 \times \bar{\Omega}^l ): \end{subarray} } \frac{1}{p+1}\frac{1}{q+1}\frac{1}{v^{(-)}+1} \sum^{p}_{i=0}\sum^{q}_{j=0}\sum^{v^{(-)}}_{k=0} \left( w(x_i,y_j,z_k,t_0)-u^\delta_{i,j,k}(t_0) \right)^2 \\
+ \varepsilon^{(-)} (t_0) \left(\left\| \frac{\partial^2 w}{\partial x^2} \right\|^2_{\begin{subarray}{c} L^2(\mathbb{R}^2 \times \bar{\Omega}^l)\end{subarray} }+ \left\| \frac{\partial^2 w}{\partial y^2} \right\|^2_{ \begin{subarray}{c} L^2(\mathbb{R}^2 \times \bar{\Omega}^l) \end{subarray}} +\left\| \frac{\partial^2 w}{\partial z^2} \right\|^2_{\begin{subarray}{c} L^2(\mathbb{R}^2 \times \bar{\Omega}^l)\end{subarray} } \right),
\end{multline}
\begin{multline}
\label{uAlphaR}
u^\varepsilon= \arg\min_{\begin{subarray}{c} w\in C^1(\mathbb{R}^2 \times \bar{\Omega}^r): \end{subarray} } \frac{1}{p+1}\frac{1}{q+1}\frac{1}{v-v^{(+)}+1} \sum^{p}_{i=0}\sum^{q}_{j=0}\sum^{v}_{k=v^{(+)}} \left( w(x_i,y_j,z_k,t_0)-u^\delta_{i,j,k}(t_0) \right)^2 \\
+ \varepsilon^{(+)} (t_0) \left(\left\| \frac{\partial^2 w}{\partial x^2} \right\|^2_{\begin{subarray}{c} L^2(\mathbb{R}^2 \times \bar{\Omega}^r) \end{subarray}}+ \left\| \frac{\partial^2 w}{\partial y^2} \right\|^2_{\begin{subarray}{c} L^2(\mathbb{R}^2 \times \bar{\Omega}^r) \end{subarray}}+\left\| \frac{\partial^2 w}{\partial z^2} \right\|^2_{\begin{subarray}{c} L^2(\mathbb{R}^2 \times \bar{\Omega}^r) \end{subarray}} \right),
\end{multline}
where the regularization parameters $\varepsilon^{(\mp)}(t_0)$ satisfy
$$ \displaystyle \frac{1}{p+1}\frac{1}{q+1}\frac{1}{v^{(-)}+1} \sum^{p}_{i=0} \sum^{q}_{j=0}\sum^{v^{(-)}}_{k=0} \left( u^\varepsilon(x_i,y_j,z_k,t_0)-u^\delta_{i,j,k}(t_0) \right)^2 = \delta^4,$$
$$ \displaystyle \frac{1}{p+1}\frac{1}{q+1}\frac{1}{v-v^{(+)}+1} \sum^{p}_{i=0}\sum^{q}_{j=0}\sum^{v}_{k=v^{(+)}} \left( u^\varepsilon(x_i,y_j,z_k,t_0)-u^\delta_{i,j,k}(t_0) \right)^2 = \delta^4.$$

Similar to \cite[Theorem 2.1]{Scherzer2018} and \cite[Theorem 3.3]{WANG2005121}, it is not difficult to derive the following assertion.

\begin{proposition}
\label{NoisyErr}
Suppose that for a.e.\ $t\in\bar{\mathcal{T}}$, $u(\cdot,\cdot,\cdot,t) \in H^2(\mathbb{R}^2 \times \Omega)$. Let $u^\varepsilon(x,y,z,t)$ be the minimizer of problems \eqref{uAlphaL} and \eqref{uAlphaR}, with $t_0$ replaced with $t$. Then, for a.e.\ $t\in\bar{\mathcal{T}}$ and $\varepsilon = \delta^2$,
\begin{align}
\label{NoisyErrIneq}
\|u^\varepsilon(\cdot,\cdot,\cdot,t) - u(\cdot,\cdot,\cdot,t)\|_{H^1(\mathbb{R}^2 \times \Omega)} \leq C_2 d^{1/4} + C_3 \sqrt{\delta}, \nonumber
\end{align}
where $d =\max \{ d_1,d_2,d_3 \}$, and $C_2$ and $C_3$ are constants that depend on the domain $\mathbb{R}^2 \times \Omega$ and the quantity $\left\| \Delta u(x,y,z,t)\right\|_{L^2(\mathbb{R}^2 \times \Omega)}$.
\end{proposition}

Using Propositions~\ref{ProAsympErr} and \ref{NoisyErr} and the triangle inequality $\|f^\delta - $ $f^*\|_{L^2(\mathbb{R}^2 \times \Omega)} \leq \|f^\delta - f_0\|_{L^2(\mathbb{R}^2 \times \Omega)} $ $ + \|f_0 - f^*\|_{L^2(\mathbb{R}^2 \times \Omega)}$, it is relatively straightforward to prove the following theorem.

\begin{theorem}
\label{ErrSource}
Defined in \eqref{fdelta}, $f^\delta$ is a stable approximation of the exact source function $f^*$ for problem (\textbf{IP}). Moreover, it has the convergence rate
\begin{equation}
\label{ErrSourceIneq}
\|f^\delta - f^*\|_{L^2(\mathbb{R}^2 \times \Omega)} = \mathcal{O} ( \mu \lvert \ln \mu \rvert + d^{1/4} + \sqrt{\delta} ).
\end{equation}
In addition, if $\mu=\mathcal{O}(\delta^{\epsilon + 1/2})$ ($\epsilon$ is any positive number) and $d^{1/4}=\mathcal{O}(\sqrt{\delta})$, then the following estimate holds:
\begin{equation*}
\|f^* - f^\delta\|_{L^2(\mathbb{R}^2 \times \Omega)} = \mathcal{O}(\sqrt{\delta}) .
\end{equation*}
\end{theorem}

Based on the above analysis, we build an efficient regularization algorithm for the three-dimensional nonlinear inverse source problem (\textbf{IP}), as shown in Algorithm~\ref{alg:Framwork}.

\begin{algorithm}[htb]
\caption{Asymptotic-expansion regularization algorithm for (\textbf{IP}).}
\label{alg:Framwork}
\begin{algorithmic}[1]
\If{The data for the full measurements $\{u^\delta_{i,j,k}, ({u_x})^\delta_{i,j,k}, ({u_y})^\delta_{i,j,k}, ({u_z})^\delta_{i,j,k}\}$ are given}
\State{break;}
\Else
\If{ The measurements $\{ u^\delta_{i,j,k} \}$ are provided}
\State Construct the smoothed data $\{u^\varepsilon,\frac{\partial u^\varepsilon}{\partial x},\frac{\partial u^\varepsilon}{\partial y},\frac{\partial u^\varepsilon}{\partial z} \}$ by solving \eqref{uAlphaL} and \eqref{uAlphaR};
\EndIf
\EndIf
\State{Calculate the approximate source function $f^\delta$ using \eqref{fdelta}.}
\end{algorithmic}
\end{algorithm}

\section{Derivation and proofs of main results} \label{derivationAndProofs}

\subsection{Proof of Theorem \ref{MainThm}}

To prove Theorem \ref{MainThm} and estimate its accuracy (\eqref{eq002} – \eqref{0orderEstim2}), we use the asymptotic method of inequalities \cite{b7}. According to this method, a solution to \eqref{mainproblem} exists if there exist continuous functions $\alpha(x, y,z, t, \mu)$ and $\beta(x, y,z, t, \mu)$, called, respectively, lower and upper solutions of \eqref{mainproblem}. First, we recall the definitions of upper and lower solutions and their role in the construction of solution  \eqref{mainproblem} \cite{b7,b8,NEFEDOV201390}.

\begin{definition} \label{Lemma1}
 The functions  $   \beta (x,y,z, t,\mu) $ and  $ \alpha (x,y,z,t, \mu)$ are called upper and lower solutions of the problem \eqref{mainproblem}, if they are continuous, twice continuously differentiable in $x$, $y$ and $z$, continuously differentiable in $t$, and, for sufficiently small $\mu$, satisfy the following conditions:

  \begin{itemize}[leftmargin=1cm]
 \item[\textbf{(C1):}] $\alpha(x,y,z,t,\mu)\leq \beta(x,y,z,t,\mu)$ for $(x,y,z,t)\in \mathbb{R}^2 \times \bar{\Omega}\times \bar{\mathcal{T}}.$
  \item[\textbf{(C2):}]  $  \displaystyle L[ \alpha]:=\mu \Delta \alpha -\frac{\partial\alpha}{\partial t}+\alpha \left( \frac{\partial\alpha}{\partial x}+\frac{\partial\alpha}{\partial y}+\frac{\partial\alpha}{\partial z}\right)-f(x,y,z)\geq 0, \quad (x,y,z,t)\in \mathbb{R}^2 \times \bar{\Omega}\times \bar{\mathcal{T}};$
  \item[\qquad \ \ ] $ \displaystyle L[\beta]:=\mu\Delta\beta-\frac{\partial\beta}{\partial t}+\beta \left( \frac{\partial\beta}{\partial x}+\frac{\partial\beta}{\partial y}+\frac{\partial\beta}{\partial z}\right)-f(x,y,z )\leq 0, \quad (x,y,z,t)\in \mathbb{R}^2 \times \bar{\Omega}\times \bar{\mathcal{T}}.$
 \item[\textbf{(C3):}] $\alpha (x,y,-a,t,\mu)\leq u^{-a} (x,y)\leq \beta(x,y,-a,t,\mu),$
 \item[\qquad \ \ ]  $ \alpha(x,y,a,t,\mu)\leq u^{a} (x,y)\leq \beta(x,y,a,t,\mu)$.
\end{itemize}
\end{definition}

%

\begin{lemma} \label{Lemma2}
(\cite{b8}) Let there be an upper $   \beta (x,y,z,t,\mu) $ and a lower $ \alpha (x,y,z,t, \mu)$ solutions to problem \eqref{mainproblem} satisfying conditions (C1)–(C3) in Definition \ref{Lemma1}. Then, under Assumptions \ref{A1}–\ref{A4}, there exists a solution $u(x,y,z,t,\mu)$ to problem \eqref{mainproblem} that satisfies the inequalities
$$
\alpha(x,y,z,t,\mu)\leq u(x,y,z,t,\mu)\leq \beta(x,y,z,t,\mu),\ (x,y,z,t)\in \mathbb{R}^2 \times \bar{\Omega}\times \bar{\mathcal{T}}.
$$
Moreover, the functions $\beta(x,y,z, t, \mu)$ and $\alpha (x,y, z,t, \mu)$ satisfy the following estimates:
\begin{align}
& \beta(x,y,z,t,\mu)-\alpha(x,y,z,t,\mu) =\mathcal{O}(\mu^{n}), \label{BetaMinusAlpha} \\
& u(x,y,z,t,\mu)=\alpha(x,y,z,t,\mu)+\mathcal{O}(\mu^{n})=U_{n-1}(x,y,z,t,\mu)+\mathcal{O}(\mu^{n}). \label{UMinusAlpha}
\end{align}
 \end{lemma}

\begin{lemma} \label{Lemma3}
(\cite{b7,NEFEDOV201390}) Lemma \ref{Lemma2} also remains valid in the case in which the functions $\alpha(x,y,z,t,\mu)$ and $\beta(x,y,z,t,\mu)$ are continuous and their derivatives with respect to $x,y,z$ have discontinuities  from the class $C^2$  in the direction of the normal to surfaces on which these solutions are not smooth, and the limit values of the derivatives on the surface $h(x,y,t)$ satisfy the following condition:
\begin{itemize}[leftmargin=1cm]
\item[\textbf{(C4):}] $\displaystyle \frac{\partial\alpha}{\partial n}(x,y,h_{\alpha}(x,y,t)+0,t,\mu)-\frac{\partial\alpha}{\partial n}(x,y,h_{\alpha}(x,y,t)-0,t,\mu)\geq 0,$
\end{itemize}
where $h_{\alpha}(x,y,t)$ is the surface on which the lower solution is not smooth;
\begin{itemize}
 \item[\qquad \   ]  $\displaystyle \frac{\partial\beta}{\partial n}(x,y,h_{\beta}(x,y,t)-0,t,\mu)-\frac{\partial\beta}{\partial n}(x,y,h_{\beta}(x,y,t)+0,t,\mu)\geq 0,$
 \end{itemize}
where $h_{\beta}(x,y,t)$ is the surface on which the upper solution is not smooth.
\end{lemma}

The proofs of Lemmas \ref{Lemma2}–\ref{Lemma3} can be found in \cite{b7,b8}. Thus, to prove Theorem \ref{MainThm}, it is necessary to construct the lower and upper solutions  $ \alpha (x,y,z,t, \mu)$  and $   \beta (x,y,z,t,\mu) $. Under conditions (C1)–(C4) for  $ \alpha (x,y,z,t, \mu)$  and $   \beta (x,y,z,t,\mu) $, estimates \eqref{NorderEstim1} and \eqref{NorderEstim2} will follow directly from Lemma \ref{Lemma2}.


We now begin the proof of Theorem \ref{MainThm}.

\begin{proof}
Following the main idea in \cite{b6},
we construct the upper and lower solutions $\beta^{(-)}$, $\beta^{(+)}$, $\alpha^{(-)}$, $\alpha^{(+)}$  and curves $h_{ \beta}$,  $h_{ \alpha}$ as a modification of asymptotic representation \eqref{asymptoticnorder}.

We introduce a positive function $ \rho (x,y,t) $, which will be defined later in \eqref{phoequat}, and use the notations  $ \rho_{\beta} (x,y,t)=-\rho (x,y,t) $ and $ \rho_{\alpha} (x,y,t)=\rho (x,y,t) $ to aid in defining the surfaces $h_{ \beta}(x,y, t)$ and $h_{ \alpha}(x,y, t)$, which in turn will determine the position of the inner transition layer for the upper and lower solution, in the form
\begin{align}\label{curveh}
\begin{dcases}
\displaystyle h_{ \beta}(x,y, t)= \sum_{i=0}^{n+1} \mu^i h_{i}(x,y,t)+\mu^{n+1} \rho_{\beta}(x,y, t),  \\
\displaystyle h_{ \alpha}(x,y, t)= \sum_{i=0}^{n+1} \mu^i h_{i}(x,y,t)+\mu^{n+1} \rho_{\alpha}(x,y, t).
\end{dcases}
\end{align}

In the vicinity of the surface $h_{\beta}(x,y, t)$,  we pass to the local coordinates $( r_{\beta},l,m)$  according to the following equations:

\begin{align*}
 &x=l-\frac{r_\beta ({h_\beta})_x}{ \sqrt{1+({h_\beta})_{x}^2+({h_\beta})_{y}^2}}, \quad y=m-\frac{r_\beta ({h_\beta})_y}{ \sqrt{1+({h_\beta})_{x}^2+({h_\beta})_{y}^2}}, \\
 &z={h_\beta}+ \frac{r_\beta  }{ \sqrt{1+({h_\beta})_{x}^2+({h_\beta})_{y}^2}}=\hat{h}_{n+1}(l,m, t)+\frac{r_\beta  }{ \sqrt{1+({h_\beta})_{x}^2+({h_\beta})_{y}^2}}+\mu^{n+1}\rho_{\beta}(l,m, t),
\end{align*}
where $r_{\beta}$ is the distance from the surface $h_{\beta}(x,y, t)$ along the normal to it, $l$ and $m$ are the coordinates of the point on the $x$ and $y$ axes respectively from which this normal is drawn, and the derivatives of the function $h_{\beta}$ at each time $t$ are taken at $x=l$ and $y=m$.

Similarly, in the vicinity of the surface $h_{\alpha}(x,y, t)$  we pass to the local coordinates \begin{align*}
 &x=l-\frac{r_\alpha ({h_\alpha})_x}{ \sqrt{1+({h_\alpha})_{x}^2+({h_\alpha})_{y}^2}}, \quad y=m-\frac{r_\alpha ({h_\alpha})_y}{ \sqrt{1+({h_\alpha})_{x}^2+({h_\alpha})_{y}^2}}, \\
 &z={h_\alpha}+ \frac{r_\alpha  }{ \sqrt{1+({h_\alpha})_{x}^2+({h_\alpha})_{y}^2}}=\hat{h}_{n+1}(l,m, t)+\frac{r_\alpha  }{ \sqrt{1+({h_\alpha})_{x}^2+({h_\alpha})_{y}^2}}+\mu^{n+1}\rho_{\alpha}(l,m, t),
\end{align*}
where $r_{\alpha}$ is the distance from the surface $h_{\alpha}(x, y, t)$ along the normal to it.

In the neighborhood of the surfaces  $h_{\beta}(x, y, t)$ and  $h_{\alpha}(x,y, t)$,  we introduce the extended variables
\begin{equation} \label{UpperLowerXi}
\xi_{\beta}=\frac{r_{\beta}}{\mu}, \quad \xi_{\alpha}=\frac{r_{\alpha}}{\mu}.
\end{equation}

The upper and lower solutions of problem \eqref{mainproblem} will be constructed separately in the domains $ \bar{D}_{ \beta}^{(-)}, \bar{D}_{ \beta}^{(+)}$ and $ \bar{D}_{ \alpha}^{(-)}, \bar{D}_{ \alpha}^{(+)}$, in which the surfaces $h_{ \beta}(x,y, t)$ and $h_{ \alpha}(x,y, t)$ divide the domain $ \mathbb{R}^2\times \bar{\Omega}\times \bar{\mathcal{T}} $ at each moment of time $t$:
\begin{align}
\label{beta}
\beta= \begin{cases}
\beta^{(-)}(x,y,z,t,\mu ), \ (x,y,z,t)\in \bar{D}_{ \beta}^{(-)}:=\{(x,y,z,t) \in \mathbb{R}^2 \times [-a, h_{\beta}] \times  \bar{\mathcal{T}} \},  \\
\beta^{(+)}(x,y,z,t,\mu ), \ (x,y,z,t)\in \bar{D}_{ \beta}^{(+)}:=\{(x,y,z,t) \in \mathbb{R}^2 \times [h_{\beta}, a] \times \bar{\mathcal{T}} \},
\end{cases}
\end{align}
\begin{align}
\label{alpha}
\alpha= \begin{cases}
\alpha^{(-)}(x,y,z,t,\mu ), \ (x,y,z,t)\in \bar{D}_{ \alpha}^{(-)}:=\{(x,y,z,t)  \in \mathbb{R}^2\times [-a, h_{ \alpha}] \times \bar{\mathcal{T}} \}, \\
\alpha^{(+)}(x,y,z,t,\mu ), \ (x,y,z,t)\in\bar{D}_{ \alpha}^{(+)}:=\{(x,y,z,t)  \in \mathbb{R}^2\times [h_{ \alpha}, a]  \times \bar{\mathcal{T}} \}.
\end{cases}
\end{align}

We match the functions $\beta^{(-)}(x,y,z, t, \mu)$, $ \beta^{(+)}(x,y, z,t, \mu)$ and $\alpha^{(-)}(x, y,z,t, \mu)$, $ \alpha^{(+)}(x,$ $ y,z,t, \mu)$ on the surfaces $h_{ \beta}(x,y, t)$ and $h_{ \alpha}(x, y,t)$, respectively, so that $\beta(x,y, t, \mu)$ and $\alpha(x, y,$ $t, \mu)$ are continuous on these surfaces and the following equations hold:
\begin{align}
\begin{dcases} \label{sewingeq}
\displaystyle \beta^{( - )}(x,y, h_{\beta}, t, \mu)=\beta^{( + )}(x,y, h_{\beta}, t, \mu)=\frac{\bar{u}^{(-)}(x,y, h_{\beta})+\bar{u}^{(+)}(x, y,h_{\beta})}{2}, \\
\displaystyle \alpha^{( - )}(x, y,h_{\alpha}, t, \mu)=\alpha^{( + )}(x,y, h_{\alpha}, t, \mu)=\frac{\bar{u}^{(-)}(x,y, h_{\alpha})+\bar{u}^{(+)}(x, y,h_{\alpha})}{2}.
\end{dcases}
\end{align}

Note that we do not match the derivatives of the upper and lower solutions on the surfaces $h_{ \beta}(x,y, t)$ and $h_{ \alpha}(x,y, t)$, respectively, so the derivatives $\partial \beta / \partial z$  and $\partial \alpha / \partial z $ have discontinuity points, and therefore fulfillment of condition (C4) is required to hold.

We construct the functions $\beta^{( \mp)}$ and $\alpha^{( \mp)}$ in the following forms:
\begin{align}\label{beta2}
\begin{dcases}
\beta^{( \mp)}= U_{n+1}^{( \mp)}|_{\xi_{\beta},h_{\beta}}+\mu^{n+1} \left(\epsilon^{( \mp)}(x,y,z)+q_{0}^{( \mp)}(\xi_{\beta},l,m,h_\beta, t)+\mu q_{1}^{( \mp)}(\xi_{\beta},l,m,h_\beta,t) \right), \\
\alpha^{( \mp)}=U_{n+1}^{( \mp)}|_{\xi_{\alpha},h_{\alpha}}-\mu^{n+1} \left( \epsilon^{( \mp)}(x,y,z)+q_{0}^{( \mp)}(\xi_{\alpha},l,m,h_\alpha,t)+\mu q_{1}^{( \mp)}(\xi_{\alpha},l,m,h_\alpha, t)\right),
\end{dcases}
\end{align}
where the functions $\epsilon^{( \mp)}(x,y,z)$ should be designed in such a way that condition (C2) is satisfied for $\beta^{( \mp)}$ and $\alpha^{( \mp)}$ in \eqref{beta2}. The functions $q_{0}^{( \mp)} $ eliminate residuals of order $\mu^n$ arising  in $L[\beta]$ and $L[\alpha]$ and the residuals of order $\mu^{n+1}$  under the condition of continuous matching of the upper solution \eqref{sewingeq}, which arise  as a result of modifying the outer part by adding  $\epsilon^{( \mp)}(x,y,z)$. The functions $q_{1}^{( \mp)} $ eliminate  residuals of order $\mu^{n+1}$  arising in $L[\beta]$ and $L[\alpha]$ with the addition of   $\epsilon^{( \mp)}(x,y,z)$ and $q_{0}^{( \mp)}  $.

We now define the functions $\epsilon^{( \mp)}(x,y,z)$  from the following equations:
\begin{align} \label{epsiloneq}
\begin{dcases}
\varphi^{( \mp)} \left(   \frac{\partial \epsilon^{( \mp)} }{\partial x} +\frac{\partial \epsilon^{( \mp)} }{\partial y}+\frac{\partial \epsilon^{( \mp)} }{\partial z}\right)   +\epsilon^{( \mp)}  \left(  \frac{\partial \varphi^{( \mp)} }{\partial x}+\frac{\partial \varphi^{( \mp)}}{\partial y}+\frac{\partial \varphi^{( \mp)}}{\partial z} \right) = -R  ,\\
\epsilon^{(-)}(x,y,-a)=R^{(-)}, \quad \epsilon^{(+)} (x,y,a)= R^{(+)},\\
\epsilon_{0}^{(\mp)}(x, y,z)=\epsilon_{0}^{(\mp)}(x+L, y,z)=\epsilon_{0}^{(\mp)}(x, y+M,z),
\end{dcases}
\end{align}
where $R, R^{(-)}, R^{(+)}$ are some positive constants, independent of $x,y,z$. The functions $\epsilon^{( \mp)} (x,$ $y,z)$ has explicit solution of the type \eqref{firstorderregularfunctionsExplicit}, where $W^{(\mp)}=\frac{-R}{\varphi^{( \mp)}}$, and since $\varphi^{(-)} (x,y,z) <0 $ and $\varphi^{(+)} (x,y,z) >0 $, we have $\epsilon^{( \mp)} (x,y,z) > 0 $ for $(x,y,z) \in \mathbb{R}^2 \times \bar{\Omega}$.

 We define the functions $q_{0}^{( \mp)}(\xi_{\beta},l,m,h_\beta,t) $ as solutions of the equation
\begin{multline} \label{q0}
\frac{\partial^2 q_{0}^{( \mp)}}{\partial \xi_{\beta}^2}+   \frac{\partial}{\partial \xi_{\beta}} \left(q_{0}^{(\mp)} \frac{{h_0}_{t}  + \tilde{u}(1-{h_0}_{l}-{h_0}_{m} ) }{ \sqrt{1+{h_0}_{l}^{2}+{h_0}_{m}^{2} }}   \right)=b_{1}^{(\mp)}(\xi_{\beta},l,m,t) \frac{\partial \rho_{\beta} }{\partial l}\\
+b_{2}^{(\mp)}(\xi_{\beta},l,m,t) \frac{\partial \rho_{\beta} }{\partial m} +b_{3}^{(\mp)}(\xi_{\beta},l,m,t) \rho_{\beta}  +b_{4}^{(\mp)}(\xi_{\beta},l,m,t)\frac{\partial \rho_{\beta} }{\partial t}+b_{5}^{(\mp)}(\xi_{\beta},l,m,t)\\
:= H_{q0}^{( \mp)}(\xi_{\beta},  l,m, t),
\end{multline}
where $b_{1}^{(\mp)}(\xi_{\beta},l,m,t)$--$b_{5}^{(\mp)}(\xi_{\beta},l,m,t)$ are known functions, and, particularly, $b_{4}^{(\mp)}(\xi_{\beta},l,m,t)=\displaystyle \frac{\partial Q_{0}^{(\mp)} (\xi_{\beta},l,m, h_0,t)}{\partial \xi_{\beta}} \frac{1}{\sqrt{1+{h_0}_{l}^{2}+{h_0}_{m}^{2}}}.$

The boundary conditions for $q_{0}^{( \mp)}(\xi_{\beta},l,m,h_\beta,t)$ follow from the equation \eqref{sewingeq}, by matching the upper solution and taking into account conditions at  $\xi_{\beta} = 0$ for the functions $Q_{i}^{( \mp)}(\xi_{\beta},l,$ $m,h_\beta, t) $ we obtain:
\begin{align}  \label{q0conditions}
\begin{dcases}
 q_{0}^{( \mp)}(0,l,m,h_\beta,t)=\frac{\epsilon^{( \pm)}(l,m,h_0)-\epsilon^{( \mp)}(l,m,h_0)}{2} := p_{q0}^{( \mp)}(l,m,h_0), \\
 q_{0}^{(\mp)}(\mp\infty,l,m,h_\beta,t) = 0.
 \end{dcases}
\end{align}

Using \eqref{q0} and \eqref{q0conditions}, we write the functions $q_{0}^{(-)}(\xi_{\beta},l,m,h_\beta, t)$ in the explicit form:
\begin{equation} \label{q0explisit}
q_{0}^{( \mp)} =J^{( \mp)}(\xi_{\beta},h_0) \left( p_{q0}^{( \mp)}(l,m,h_0)   +\int_{0}^{\xi_{\beta}} \frac{1}{J^{( \mp)}(s,h_0)} \int_{\mp \infty}^{s} H_{q0}^{( \mp)} (\eta,l,m,t) d\eta ds \right).
\end{equation}

 The functions $q_{1}^{( \mp)}(\xi_{\beta},l,m,h_\beta, t) $ are determined from the following equation:
\begin{align} \label{q1}
\begin{split}
\frac{\partial^2 q_{1}^{( \mp)}}{\partial \xi_{\beta}^2}+   \frac{\partial}{\partial \xi_{\beta}} \left(q_{1}^{(\mp)} \frac{{h_0}_{t}   + \tilde{u}(1-{h_0}_{l} -{h_0}_{m}) }{ \sqrt{1+{h_0}_{l}^{2}+{h_0}_{m}^{2} }}   \right)= H_{q1}^{( \mp)}(\xi_{\beta}, l,m, t),
\end{split}
\end{align}
where $H_{q1}^{( \mp)}(\xi_{\beta},l,m, t) $ depend on known functions $ h_{0,1}$, $u_{0,1}^{\mp} $, $Q_{0,1}^{\mp}$,  $\rho_{\beta}$, $\epsilon^{( \mp)}$, and $ q_{0}^{( \mp)}$. For equation \eqref{q1}, we infer that the boundary condition is
\begin{align*}
q_{1}^{( \mp)}(0,l,m,h_\beta,t)=0, \ q_{1}^{( \mp)}(\xi_{\beta},l,m,h_\beta,t) \rightarrow 0 \ \text{for} \  \xi_{\beta} \rightarrow \mp \infty.
\end{align*}

Replacing $\rho_{\beta}$ with $\rho_{\alpha}$ and $\xi_{\beta}$ with $\xi_{\alpha}$ in \eqref{q0}–\eqref{q1}, we define the functions $q_{0}^{( \mp)}(\xi_{\alpha},l,m,$ $h_\alpha, t)$ and $q_{1}^{( \mp)}(\xi_{\alpha},l,m,h_\alpha, t)$ that appear in the functions $\alpha^{( \mp)}$.

The functions $q_{0}^{( \mp)}$ and $q_{1}^{( \mp)}$ satisfy exponential estimates of the types in \eqref{equat22} and \eqref{equat23}.

Now, we have to show that the functions $\beta(x,y,z,t,\mu)$ and $\alpha(x,y,z,t,\mu)$ are upper and lower solutions to problem \eqref{mainproblem}. To do this, we check conditions (C1)–(C4).

First, we check that condition (C1) has been fulfilled, with regard to the ordering of the lower and upper solutions. To do this, we consider three regions that illustrate the difference between the upper and lower solutions, $\beta-\alpha $:
\begin{align}
\quad \beta-\alpha=\left\{\begin{array}{l}
\beta^{(-)}-\alpha^{(-)}, \quad \RomanNumeralCaps{1} =  \{ (x,y,z,t):   \mathbb{R}^2 \times [-a, h_{\beta}(x,y, t)] \times \bar{\mathcal{T}} \},\\
\beta^{(+)}-\alpha^{(-)} , \quad \RomanNumeralCaps{2} = \{ (x,y,z,t):  \mathbb{R}^2 \times [ h_{\beta}(x,y, t), h_{\alpha}(x,y, t)] \times \bar{\mathcal{T}} \},\\
\beta^{(+)}-\alpha^{(+)} , \quad \RomanNumeralCaps{3} = \{ (x,y,z,t):  \mathbb{R}^2 \times [h_{\alpha}(x,y, t), a] \times \bar{\mathcal{T}} \}.
\end{array}\right.
\end{align}
First, we find the variables on which $\beta$ and $\alpha $ depend.
\begin{align*}
z&=\hat{h}_{n+1}(l,m, t)+\frac{r_\beta  }{ \sqrt{1+(\hat{h}_{n+1})_{x}^2+(\hat{h}_{n+1})_{y}^2}}-\mu^{n+1}\rho(l,m, t)+\mathcal{O}(\mu^{n+1})\\
&=\hat{h}_{n+1}(l,m, t)+\frac{r_\alpha  }{ \sqrt{1+(\hat{h}_{n+1})_{x}^2+(\hat{h}_{n+1})_{y}^2}}+\mu^{n+1}\rho(l,m, t)+\mathcal{O}(\mu^{n+1}).
\end{align*}

We rewrite \eqref{UpperLowerXi} in the following form:

\begin{align*}
\displaystyle \xi_{\beta}=\frac{z-h_\beta}{\mu }  \sqrt{1+({h_\beta})_{x}^{2}+({h_\beta})_{y}^{2}}, \quad \xi_{\alpha}=\frac{z-h_\alpha}{\mu }  \sqrt{1+({h_\alpha})_{x}^{2}+({h_\alpha})_{y}^{2}},
\end{align*}
from where we find:
$$
\Delta \xi=\xi_{\beta}-\xi_{\alpha}=2\mu^{n} \left( \rho\sqrt{1+(h_{0})_{x}^{2}+(h_{0})_{y}^{2}}+\frac{\left( \rho_x ({h_0})_{x}+\rho_y ({h_0})_{y} \right)(h_0-y)}{\sqrt{1+(h_{0})_{x}^{2}+(h_{0})_{y}^{2}}}\right)+\mathcal{O}(\mu^{n+1}).
$$

For region \RomanNumeralCaps{2}, the following holds:
\begin{align*}
    \begin{split}
       & 0\leq\xi_{\beta}\leq 2\mu^{n}\rho(l,m, t)\sqrt{1+(h_{0})_{x}^{2}+(h_{0})_{y}^{2}}+\mathcal{O}(\mu^{n+1}),\\
       & 0 \geq\xi_{\alpha}\geq  -2\mu^{n}\rho(l,m, t)\sqrt{1+(h_{0})_{x}^{2}+(h_{0})_{y}^{2}}+\mathcal{O}(\mu^{n+1}),\\
       & \xi_{\beta}-\xi_{\alpha}=2\mu^{n}  \rho(l,m, t)\sqrt{1+(h_{0})_{x}^{2}+(h_{0})_{y}^{2}}+\mathcal{O}(\mu^{n+1}).
    \end{split}
\end{align*}

We can write an expression for the difference between the upper and lower solutions:
\begin{align} \label{diffbetaalpha1}
\begin{array}{ll}
\beta^{(+)}-\alpha^{(-)} =  \\
 \quad \sum_{i=0}^{n}\mu^{i} \left(\bar{u}_{i}^{(+)}\left(l-\frac{\mu\xi_{\beta} ({h_\beta})_x}{ \sqrt{1+({h_\beta})_{x}^2+({h_\beta})_{y}^2}}, m-\frac{\mu\xi_{\beta} ({h_\beta})_y}{ \sqrt{1+({h_\beta})_{x}^2+({h_\beta})_{y}^2}},{h_\beta}+ \frac{\mu\xi_{\beta}   }{ \sqrt{1+({h_\beta})_{x}^2+({h_\beta})_{y}^2}}\right) \right. \\
\quad -\left.\bar{u}_{i}^{(-)}\left(l-\frac{\mu\xi_{\alpha} ({h_\alpha})_x}{ \sqrt{1+({h_\alpha})_{x}^2+({h_\alpha})_{y}^2}}, m-\frac{\mu\xi_{\alpha} ({h_\alpha})_y}{ \sqrt{1+({h_\alpha})_{x}^2+({h_\alpha})_{y}^2}},{h_\alpha}+ \frac{\mu\xi_{\alpha}   }{ \sqrt{1+({h_\alpha})_{x}^2+({h_\alpha})_{y}^2}}\right)\right)\\
\quad +\sum_{i=0}^{n}\mu^{i} \left( Q_{i}^{(+)}\left(\xi_{\beta}, l,m, h_{\beta}, t \right)  -Q_{i}^{(-)} \left(\xi_{\alpha}, l,m, h_{\alpha}, t \right) \right) +\mathcal{O}(\mu^{n+1}).
\end{array}
\end{align}

Expanding equation \eqref{diffbetaalpha1} in series, and taking into account the notation \eqref{derivativetildeu},  equation \eqref{sewindcondexpanded0} and the fact that in region \RomanNumeralCaps{2} $\xi_{\beta}=\mathcal{O}(\mu^{n})$ and  $ \xi_{\alpha}=\mathcal{O}(\mu^{n})$, we obtain an expression for the difference between the upper and lower solutions in region \RomanNumeralCaps{2}:
\begin{align}
\begin{split}
\beta^{(+)}-\alpha^{(-)} &=\frac{\partial Q_{0}^{(+)}}{\partial \xi}(0, x,y, h_{0}(x,y, t),t)\xi_{\beta}-\frac{\partial Q_{0}^{(-)}}{\partial \xi}(0, x,y, h_{0}(x,y, t),t)\xi_{\alpha}+\mathcal{O}(\mu^{n+1})\\
&=2\mu^{n}\rho(l,m, t)\sqrt{1+(h_{0})_{x}^{2}+(h_{0})_{y}^{2}}\cdot\Phi^{(+)}(0, h_{0}(x,y, t))+\mathcal{O}(\mu^{n+1}).
\end{split}
\end{align}

 Using equation \eqref{Q0equation} and Assumption \ref{A3}, it is possible to verify that $\Phi^{(+)}(0, h_{0}(x,y, t))>0$, and, for positive values of $\rho$ and for a sufficiently small $\mu$, we obtain $\beta-\alpha>0, \ (x,y,z,t)\in \RomanNumeralCaps{2}.$

Now let us consider the difference between the upper and lower solutions at region \RomanNumeralCaps{3}, where $\xi_{\alpha}\geq 0, \ \xi_{\beta}=\xi_{\alpha}+\Delta \xi. $ Using exponential properties of the functions $Q_{i}^{(+)}$ and $q_{0}^{(+)}$, we obtain:
\begin{multline}
\displaystyle \beta-\alpha=\beta^{(+)}-\alpha^{(+)}=\sum_{i=0}^{n}\mu^{i}\left(Q_{i}^{(+)}(\xi_{\beta}, l,m, h_{\beta}, t)-Q_{i}^{(+)}(\xi_{\alpha}, l,m, h_{\alpha}, t)\right)\\
+2\mu^{n+1}\epsilon^{(+)}+\displaystyle \mu^{n+1}\left(q_{0}^{(+)}(\xi_{\beta},l,m,h_{\beta}, t)-q_{0}^{(+)}(\xi_{\alpha},l,m,h_{\alpha}, t)\right)+\mathcal{O}(\mu^{n+2})\\
=2\mu^{n+1}\epsilon^{(+)}+\frac{\partial Q_{0}^{(+)}}{\partial\xi}(\xi_{\alpha}, l,m, h_{\alpha}, t)(\xi_{\beta}-\xi_{\alpha})+\mathcal{O}(\mu^{n+1})\exp(-\kappa_{1}\xi_{\alpha})+\mathcal{O}(\mu^{n+2}),
\end{multline}
where $\kappa_{1}>0$ is a constant independent of $\xi_{\alpha},l,m,t$.

Taking into account estimates \eqref{equat22}  and \eqref{equat23} and the equality $ \xi_{\beta}-\xi_{\alpha}=\mathcal{O}(\mu^{n}) $, we obtain an expression for the difference between the upper and lower solutions in region \RomanNumeralCaps{3}:
\begin{equation} \label{diffbetaalpha2}
\beta-\alpha\leq 2\mu^{n+1}\epsilon^{(+)}+\{C_{4}\mu^{n}\exp(-\kappa_{0}\xi_{\alpha})-C_{5}\mu^{n+1}\exp(-\kappa_{1}\xi_{\alpha})\}+\mathcal{O}(\mu^{n+2}),
\end{equation}
where $\kappa_{0}>0$, $C_{4}>0$ and $C_{5}>0$ are some constants independent of $\xi_{\alpha},l,m,t$.

Let $\kappa_{0}\leq\kappa_{1}$, the expression in the brackets in \eqref{diffbetaalpha2} is positive, since $ C_{4}>C_{5}\mu$ for a sufficiently small $\mu$. Hence, $\beta-\alpha>0$.

Let $\kappa_{0}>\kappa_{1}$. Consider the region $x\in \mathbb{R},y\in \mathbb{R}, z \in [h_{\alpha},h_{\alpha}+\frac{N\mu}{ \sqrt{1+({h_\alpha})_{x}^{2}+({h_\alpha})_{y}^{2}}}], t\in \bar{\mathcal{T}} ,$ where $N>0$. In this region, the value   $r_{\alpha}$ changes on the interval $[0, N\mu]$ and the inequality $\exp(-\kappa_{0}\xi_{\alpha})\geq\exp(-\kappa_{0}N)$ is satisfied, and so the bracketed expression in \eqref{diffbetaalpha2} is positive if $\mu$ is small enough due to the component $C_{4}\mu^{n}\exp(-\kappa_{0}\xi_{\alpha})$. Hence,  $\beta-\alpha>0.$

Now we choose a number $N$ large enough to satisfy the inequality $C_{5}\exp(-\kappa_{1}  N)<2\epsilon^{(+)}$. When $ z \in [h_{\alpha}+\frac{N\mu}{ \sqrt{1+({h_\alpha})_{x}^{2}+({h_\alpha})_{y}^{2}}}, a]  $, due to the choice of the number $N$ we obtain
$$
2\mu^{n+1}\epsilon^{(+)}-C_{5}\mu^{n+1}\exp(-\kappa_{1}\xi_{\alpha})\geq\mu^{n+1}(2\epsilon^{(+)}-C_{5}\exp(-\kappa_{1}N))>0.
$$

Thus, $\beta(x, y, z,t, \mu)-\alpha(x, y,z, t, \mu)>0$ everywhere in region \RomanNumeralCaps{3}. The proof of the inequality $\beta(x, y,z, t, \mu)-\alpha(x, y,z, t, \mu)>0$ for region \RomanNumeralCaps{1} is developed in the same way as for region \RomanNumeralCaps{3}.

The method of constructing the upper and lower solutions implies the inequalities
\begin{equation*}
L[\beta]=-\mu^{n+1} R + \mathcal{O}(\mu^{n+2})<0, \quad L[\alpha]=\mu^{n+1} R + \mathcal{O}(\mu^{n+2})>0,
\end{equation*}
where $R$ is a constant from \eqref{epsiloneq}. This verifies condition (C2).

Condition (C3) is satisfied for sufficiently large values $R^{(-)}$ and $R^{(+)}$ in the boundary conditions of equation \eqref{epsiloneq}.

We now check condition (C4) for the upper solution, and expand it in powers of $\mu$; due to the matching of the asymptotics \eqref{sewindcondexpanded0} and \eqref{matchingfirstord}, the coefficients at $\mu^{i}$ for $i = 1,\ldots,n$ are equal to zero, and the coefficient at $\mu^{n + 1}$ includes only terms that arise as a result of the modification of the asymptotic expansion:
\begin{align}
\mu \left( \frac{\partial\beta^{( - )}}{\partial n}-\frac{\partial\beta^{( + )}}{\partial n} \right)\Big|_{h=h_{\beta}} = \mu^{n+1} \left( \frac{\partial {q_0}^{( - )}}{\partial \xi_{\beta}}(0,x,y,h_\beta,t) - \frac{\partial q_{0}^{( + )}}{\partial \xi_{\beta}}(0,x,y,h_\beta,t) \right) +\mathcal{O}(\mu^{n+2}).
\end{align}

Using the explicit solution for $  {q_0}^{( \mp)}(0,x,y,h_\beta,t)$ \eqref{q0explisit}, we find
\begin{multline} \label{q0lminusq0r}
 \frac{\partial {q_0}^{( - )}}{\partial \xi_{\beta}}(0,x,y,h_\beta,t) - \frac{\partial q_{0}^{( + )}}{\partial \xi_{\beta}}(0,x,y,h_\beta,t) \\
 =\frac{\partial \rho_{\beta}(x,y, t) }{\partial t} \left( \int_{- \infty}^{0} b_{4}^{(-)}(s,x,y,t) ds \right.\left. + \int_{0}^{+ \infty}b_{4}^{(+)}(s,x,y,t) ds \right) +\frac{\partial \rho_{\beta}(x,y, t) }{\partial l} H_1(x,y,t)
\\   +\frac{\partial \rho_{\beta}(x,y, t) }{\partial m} H_2(x,y,t)+\rho_{\beta} (x,y, t) H_3(x,y,t)+H_4(x,y,t),
\end{multline}
where
\begin{flalign*}
H_1(x,y,t)=  \int_{- \infty}^{0} b_{1}^{(-)}(  s,x,y,t) ds + \int_{0}^{+ \infty} b_{1}^{(+)}(s,x,y,t) ds  ,\\
H_2(x,y,t)=  \int_{- \infty}^{0} b_{2}^{(-)}(s,x,y,t) ds + \int_{0}^{+ \infty} b_{2}^{(+)}(s,x,y,t) ds ,\\
H_3(x,y,t)=  \int_{- \infty}^{0} b_{3}^{(-)}(s,x,y,t) ds + \int_{0}^{+ \infty} b_{3}^{(+)}(s,x,y,t) ds ,
\end{flalign*}
\begin{align*}
H_4(x,y,t)=\left(\epsilon^{( +)} -\epsilon^{( -)} \right) \left(  \frac{- {h_0}_{t}   + \phi_0 ({h_0}_{x}+{h_0}_{y} -1) }{ \sqrt{1+{h_0}_{x}^{2}+{h_0}_{y}^{2} }} \right) \\
+\int_{- \infty}^{0} b_{5}^{(-)}(s,x,y,t) ds + \int_{0}^{+ \infty} b_{5}^{(+)}(s,x,y,t) ds.
\end{align*}

The function $\rho_{\beta}(x,y, t) $ is obtained as a solution to the problem
\begin{align} \label{phoequat}
\begin{dcases}
\frac{\partial \rho_{\beta} }{\partial t} \frac{\varphi^{(+)}-\varphi^{(-)}}{\sqrt{1+{h_0}_{x}^{2}+{h_0}_{y}^{2}}}  
 = -H_1(x,y,t) \frac{\partial \rho_{\beta} }{\partial x}-H_2(x,y,t) \frac{\partial \rho_{\beta} }{\partial y}   -H_3(x,y,t)\rho_{\beta}-H_4(x,y,t)+\sigma, \\
\rho_{\beta} (x,y,0)=\rho^0(x,y), \quad  \rho_{\beta} (x,y,t)=\rho_{\beta} (x+M,t)=\rho_{\beta} (x,y+L,t), \quad t \in \bar{\mathcal{T}}.
\end{dcases}
\end{align}
Since the difference $\varphi^{(+)}-\varphi^{(-)} $ is positive, and we choose the constant  $\sigma$ and the function $\rho^0(x,y) $  to be positive for any $x,y$, then   the solution $\rho (x,y,t)$ of the equation \eqref{phoequat} is also positive for a sufficiently large $\sigma$.

For such $\rho (x,y,t)$, we obtain:
\begin{equation}
\mu \left( \frac{\partial\beta^{( - )}}{\partial n}-\frac{\partial\beta^{( + )}}{\partial n} \right)\Big|_{h(l,m,t)=h_{\beta}(l,m,t)} =\mu^{n+1} \sigma +\mathcal{O}(\mu^{n+2}) >0.
\end{equation}

Similarly, condition (C4) is satisfied for the functions $\alpha^{( \mp)}$, and the constructed upper and lower solutions guarantee the existence of a solution $u (x,y,z, t, \mu)$ to the problem \eqref{mainproblem} satisfying the inequalities
\begin{equation}
 \alpha(x,y,z,t,\mu)\leq u(x,y,z,t, \mu)\leq \beta(x,y,z,t,\mu),
\end{equation}
by addition, estimates \eqref{NorderEstim1} and \eqref{NorderEstim2} are valid.

It is possible to show that the estimate \eqref{NorderEstim3} is also holds using technique represented in \cite{NEFEDOV201390} by estimating the difference  $Z_{n}(x,y, t,$ $ \mu)\equiv u(x,y, t, \mu)-U_{n}(x,y, t, \mu) $;
 the function $Z_{n}(x,y,t, \mu)$ satisfies the equation
 \begin{align} \label{proof2.22}
\Delta Z_n- \frac{1}{\mu} \frac{\partial z_{n}}{\partial t}-K Z_{n}  =-K Z_{n}+ \frac{1}{\mu}\left(\frac{\partial}{\partial x}+\frac{\partial}{\partial y}+\frac{\partial}{\partial z} \right)\int_{u}^{U_{n}} ( s) ds + \mu^{n}\psi(x,y, t, \mu),
\end{align}
for $(x,y, t)\in \mathbb{R} \times \bar{\Omega}\times \bar{\mathcal{T}}$, with zero boundary conditions, where $|\psi(x,y, t, \mu)|\leq c_{0}$ and $ c_{0}$  is a constant independent of $x,y, z, t, \mu $. Using the estimates from Lemma \ref{Lemma2}, and solving \eqref{proof2.22} with a Green's function for the parabolic operator  \cite{pao1992}, taking into account the estimates \cite[Page 49]{pao1992} we obtain $\displaystyle \frac{\partial Z_{n}}{\partial x}=\mathcal{O}(\mu^{n})$, $\displaystyle \frac{\partial Z_{n}}{\partial y}=\mathcal{O}(\mu^{n})$,$\displaystyle \frac{\partial Z_{n}}{\partial z}=\mathcal{O}(\mu^{n})$ for $(x,y, t)\in \mathbb{R} \times \bar{\Omega}\times \bar{\mathcal{T}}$ and the estimate $\displaystyle \frac{\partial Z_{n}}{\partial n}=\mathcal{O}(\mu^{n})$ follows directly from \eqref{derivativetonormalXYZ}. This completes the proof of Theorem \ref{MainThm}.
\end{proof}

\subsection{Proof of Proposition 1}

\begin{proof}
First, we note that the exact source function $f^*$ has the following representation according to equations \eqref{zeroorderregularequation1} and \eqref{zeroorderregularequation2}:

\begin{align} \label{equatforf*}
f^*= \begin{cases}
\displaystyle \bar{u}_{0}^{(-)} \left(\frac{\partial\bar{u}_{0}^{(-)}}{\partial x}+\frac{\partial\bar{u}_{0}^{(-)}}{\partial y}+\frac{\partial\bar{u}_{0}^{(-)}}{\partial z}\right),  \quad (x,y,z)\in \mathbb{R}^2 \times  (-a,h_0(x,y,t) - \Delta h/2), \\ \\
\displaystyle \bar{u}_{0}^{(+)}\left( \frac{\partial\bar{u}_{0}^{(+)}}{\partial x}+\frac{\partial\bar{u}_{0}^{(+)}}{\partial y}+\frac{\partial\bar{u}_{0}^{(-)}}{\partial z}\right), \quad (x,y,z)\in \mathbb{R}^2 \times  (h_0(x,y,t) + \Delta h/2,a).
  \end{cases}
\end{align}

Let $\Omega'=(-a,h_0(x,y,t) - \Delta h/2)$. By using Corollary \ref{Corollary1}, we deduce, together with $\lvert \Omega' \rvert \leq \lvert \Omega \rvert=a$, that

\begin{align} \label{eq004IPPf}
\begin{dcases}
&\left\| \varphi^{(-)}(x,y,z) - u (x,y,z,t) \right\|_{L^p(\mathbb{R}^2 \times\Omega')} \leq C \mu, \quad \left\| \frac{\partial \varphi^{(-)}}{\partial x} - \frac{\partial  u}{\partial x} \right\|_{L^p(\mathbb{R}^2 \times\Omega')} \leq C \mu , \\
& \left\| \frac{\partial  \varphi^{(-)}}{\partial y} - \frac{\partial  u}{\partial y} \right\|_{L^p(\mathbb{R}^2 \times\Omega')} \leq C \mu, \quad \left\| \frac{\partial  \varphi^{(-)}}{\partial z} - \frac{\partial  u}{\partial z} \right\|_{L^p(\mathbb{R}^2 \times\Omega')} \leq C \mu,\\
&\left\| \frac{\partial  \varphi^{(-)}}{\partial x} \right\|_{L^p(\mathbb{R}^2 \times \Omega')} \leq \left\| \frac{\partial  \varphi^{(-)}}{\partial x} - \frac{\partial  u}{\partial x} \right\|_{L^p(\mathbb{R}^2 \times\Omega')} + \left\|\frac{\partial  u}{\partial x} \right\|_{L^p(\mathbb{R}^2 \times \Omega')}\leq C \mu + \left\|\frac{\partial  u}{\partial x} \right\|_{L^p(\mathbb{R}^2 \times \Omega')},\\
&\left\| \frac{\partial  \varphi^{(-)}}{\partial y} \right\|_{L^p(\mathbb{R}^2 \times \Omega')} \leq \left\| \frac{\partial  \varphi^{(-)}}{\partial y} - \frac{\partial  u}{\partial y} \right\|_{L^p(\mathbb{R}^2 \times \Omega')} + \left\|\frac{\partial  u}{\partial y} \right\|_{L^p(\mathbb{R}^2 \times \Omega')} \leq C \mu + \left\|\frac{\partial  u}{\partial y} \right\|_{L^p(\mathbb{R}^2 \times \Omega')},\\
&\left\| \frac{\partial  \varphi^{(-)}}{\partial z} \right\|_{L^p(\mathbb{R}^2 \times \Omega')} \leq \left\| \frac{\partial  \varphi^{(-)}}{\partial z} - \frac{\partial  u}{\partial z} \right\|_{L^p(\mathbb{R}^2 \times \Omega')} + \left\|\frac{\partial  u}{\partial z} \right\|_{L^p(\mathbb{R}^2 \times \Omega')} \leq C \mu + \left\|\frac{\partial  u}{\partial z} \right\|_{L^p(\mathbb{R}^2 \times \Omega')}.
\end{dcases}
\end{align}

From estimates \eqref{eq004IPPf}, we conclude that
\begin{align} \label{proofPro1}
\begin{split}
& \left\|f^* - f_0 \right\|_{L^p(\mathbb{R}^2 \times \Omega')}  = \left\| \varphi^{(-)} \left(  \frac{\partial \varphi^{(-)}}{\partial x} +\frac{\partial \varphi^{(-)}}{\partial y} +\frac{\partial \varphi^{(-)}}{\partial z} \right) - u  \left(  \frac{\partial u }{\partial x} + \frac{\partial u }{\partial y} + \frac{\partial u }{\partial z} \right) \right\|_{L^p(\mathbb{R}^2 \times \Omega')} \\
& \quad \leq \left\| \varphi^{(-)} \left(  \frac{\partial \varphi^{(-)}}{\partial x} +\frac{\partial \varphi^{(-)}}{\partial y}+\frac{\partial \varphi^{(-)}}{\partial z} \right) - u  \left(  \frac{\partial \varphi^{(-)} }{\partial x} + \frac{\partial \varphi^{(-)} }{\partial y} + \frac{\partial \varphi^{(-)} }{\partial z}\right) \right\|_{L^p(\mathbb{R}^2 \times \Omega')}   \\
& \qquad \qquad \qquad   +\left\| u \left(  \frac{\partial \varphi^{(-)}}{\partial x} +\frac{\partial \varphi^{(-)}}{\partial y}+\frac{\partial \varphi^{(-)}}{\partial z} \right) - u  \left(  \frac{\partial u }{\partial x} + \frac{\partial u }{\partial y}+ \frac{\partial u }{\partial z} \right) \right\|_{L^p(\mathbb{R}^2 \times \Omega')} \\
 & \quad \leq \left\| \varphi^{(-)} - u  \right\|_{L^p(\mathbb{R}^2 \times \Omega')} \left\|  \frac{\partial \varphi^{(-)}}{\partial x} +\frac{\partial \varphi^{(-)}}{\partial y}  +\frac{\partial \varphi^{(-)}}{\partial z} \right\|_{L^p(\mathbb{R}^2 \times \Omega')} \\
& \qquad \qquad \qquad    +   \left\| u  \right\|_{L^p(\mathbb{R}^2 \times \Omega')} \left\|   \frac{\partial \varphi^{(-)}}{\partial x} +\frac{\partial \varphi^{(-)}}{\partial y}+\frac{\partial \varphi^{(-)}}{\partial z} - \frac{\partial u }{\partial x} - \frac{\partial u }{\partial y}- \frac{\partial u }{\partial z} \right\|_{L^p(\mathbb{R}^2 \times \Omega')}\\
 & \quad   \leq 3C  \left( C \mu +\left\|u \right\|_{W^{1,p}(\mathbb{R}^2 \times \Omega')} \right) \mu =: c_1  \mu.
\end{split}
\end{align}

Following exactly the same reasoning process, we also derive the following inequality:
\begin{equation} \label{proofPro2}
\left\|f^* - f_0 \right\|_{L^p(\mathbb{R}^2 \times [h_0 + \Delta h/2,a])}\leq c_2 \mu
\end{equation}
with a constant $c_2$. Now, we note that
\begin{multline} \label{proofPro3}
\displaystyle \left\|f^* - f_0 \right\|_{L^p(\mathbb{R} \times [h_0 - \Delta h/2, h_0 + \Delta/2])}\\
\leq
\left\|f^* \right\|_{L^p(\mathbb{R}^2 \times [h_0 - \Delta h/2, h_0 + \Delta h/2])}
+ \left\| f_0 \right\|_{L^p(\mathbb{R}^2 \times [h_0 - \Delta h/2, h_0 + \Delta h/2])} \leq c_3 \mu \lvert \ln \mu \rvert
\end{multline}
with $c_3=\left\|f^* \right\|_{C(\mathbb{R}^2 \times \Omega)} + \left\| f_0 \right\|_{C(\mathbb{R}^2 \times \Omega)}$. By combining \eqref{proofPro1}–\eqref{proofPro3}, we deduce that
\begin{multline*}
\left\|f^* - f_0 \right\|^p_{L^p(\mathbb{R}^2 \times [-a,a])}  =  \left\|f^* - f_0 \right\|^p_{L^p(\mathbb{R}^2 \times [-a,h_0 - \frac{\Delta h}{2}])} + \left\|f^* - f_0 \right\|^p_{L^p(\mathbb{R}^2 \times [h_0 + \frac{\Delta h}{2},a] )} \\
 + \left\|f^* - f_0 \right\|^p_{L^p(\mathbb{R}^2 \times [h_0 - \frac{\Delta h}{2},h_0 + \frac{\Delta h}{2}])}   \leq c^p_1 \mu^p + c^p_2 \mu^p + c^p_3 \mu^p \lvert \ln \mu \rvert^p,
\end{multline*}
which yields the required estimate \eqref{f0Ineq}, with $C_1 = \left( c^p_1  + c^p_2 + c^p_3 \right)^{1/p}$.
\end{proof}

\subsection{Proof of Theorem \ref{ErrSource}}

\begin{proof}
Let $f_0$ be the pre-approximate source function defined in \eqref{f0}. Similarly to the proof of Proposition \ref{ProAsympErr}, for small enough $\delta$ we obtain the following inequalities:
\begin{align}
\label{f0falphaErr}
\begin{array}{ll}
\displaystyle  \|f^\delta - f_0  \|_{L^2(\mathbb{R}^2 \times \Omega')}  = \left\| u^\delta \left(  \frac{\partial u^\delta }{\partial x} +\frac{\partial u^\delta }{\partial y}+\frac{\partial u^\delta }{\partial z} \right) - u  \left(  \frac{\partial u }{\partial x} + \frac{\partial u }{\partial y}+ \frac{\partial u }{\partial z} \right) \right\|_{L^2(\mathbb{R}^2 \times \Omega')}   \\
 \displaystyle  \ \leq \left\| u^\delta - u  \right\|_{L^2(\mathbb{R}^2 \times \Omega')} \left\|  \frac{\partial u^\delta}{\partial x} +\frac{\partial u^\delta}{\partial y}+\frac{\partial u^\delta}{\partial z} \right\|_{L^2(\mathbb{R}^2 \times \Omega')}   \\
  \displaystyle \qquad \qquad \qquad \qquad  \qquad \qquad  +   \left\| u  \right\|_{L^2(\mathbb{R}^2 \times \Omega')} \left\|   \frac{\partial u^\delta}{\partial x} +\frac{\partial u^\delta}{\partial y}+\frac{\partial u^\delta}{\partial z} - \frac{\partial u }{\partial x} - \frac{\partial u }{\partial y}- \frac{\partial u }{\partial z} \right\|_{L^2(\mathbb{R}^2 \times \Omega')}\\
 \displaystyle  \  \leq \left\| u^\delta - u  \right\|_{L^2(\mathbb{R}^2 \times \Omega')}  \left( \|u^\delta - u   \|_{H^{1}(\mathbb{R}^2 \times \Omega')} + \|u   \|_{H^{1}(\mathbb{R}^2 \times \Omega')} \right)  +   \left\| u  \right\|_{L^2(\mathbb{R}^2 \times \Omega')}   \|u^\delta - u   \|_{H^{1}(\mathbb{R}^2 \times \Omega')} \\
\displaystyle  \
\leq 2 \left( \|  u \|_{H^{1}(\mathbb{R}^2 \times \Omega')} + 1 \right)  \|u^\delta - u   \|_{H^{1}(\mathbb{R}^2 \times \Omega')} \\  \displaystyle  \
\leq 2 \left( \|  u \|_{H^{1}(\mathbb{R}^2 \times \Omega')} + 1 \right) \left( C_2 d^{1/4} + C_3 \sqrt{\delta}  \right)  \\ \displaystyle  \
\leq C_6 (d^{1/4} + \sqrt{\delta}),
\end{array}
\end{align}
where $\displaystyle  C_6 = 2 \left( \|  u \|_{H^{1}(\mathbb{R}^2 \times \Omega')} + 1 \right) \left(   C_2  + C_3  \right)$.

By combining \eqref{f0Ineq} and \eqref{f0falphaErr}, we can derive the estimate \eqref{ErrSourceIneq} from the triangle inequality $\|f^\delta - f^*  \|_{L^2(\mathbb{R}^2 \times\Omega)}\leq \|f^\delta - f_0 \|_{L^2(\mathbb{R}^2 \times\Omega)} + \|f_0 -f^* \|_{L^2(\mathbb{R}^2 \times\Omega)}  $.
\end{proof}

\section{Numerical example} \label{simulation}

\subsection{Forward problem}

We consider the following three-dimensional RDA equation with the given source function $f(x,y,z)=\cos{ ( \pi x )} \cos{( \pi y )} \cos{( \pi z )}$:
\begin{align} \label{forwardexample1}
\begin{dcases}
\displaystyle 0.01 \Delta u - \frac{\partial u}{\partial t} = -u \left( \frac{\partial u}{\partial x} + \frac{\partial u}{\partial y}+ \frac{\partial u}{\partial z} \right) +f(x,y,z), \\
\displaystyle u(x,y,-1,t)=-4, \quad u(x,y,1,t)=2, \quad u(x,y,z,t)=u(x+2,y,z,t)=u(x,y+2,z,t), \\
\displaystyle u(x,y,z,0)=u_{init}(x,y,z,\mu ), \quad \\
x \in [-1, 1], y \in [-1, 1], z \in [-1, 1], t \in [0,0.5].
\end{dcases}
\end{align}

Using the asymptotic method, we reduce the original equation to two problems that determine the outer functions:
\begin{align}\label{eqUL}
\begin{cases}
\displaystyle \varphi^{(-)} \left( \frac{\partial \varphi^{(-)}}{\partial x} + \frac{\partial \varphi^{(-)} }{\partial y}+ \frac{\partial \varphi^{(-)} }{\partial z} \right) = \cos{ ( \pi x )} \cos{( \pi y )} \cos{( \pi z )}, \\
\varphi^{(-)}(x,y,-1)=-4, \ \varphi^{(-)}(x,y,z)=\varphi^{(-)}(x+2,y,z)=\varphi^{(-)}(x,y+2,z),
\end{cases}
\end{align}
\begin{align}\label{eqUR}
\begin{cases}
\displaystyle \varphi^{(+)} \left( \frac{\partial \varphi^{(+)} }{\partial x} + \frac{\partial \varphi^{(+)} }{\partial y}+ \frac{\partial \varphi^{(+)} }{\partial z} \right) = \cos{ ( \pi x )} \cos{( \pi y )} \cos{( \pi z )} , \\
\displaystyle \varphi^{(+)}(x,y,1)=2, \ \varphi^{(+)}(x,y,z)=\varphi^{(+)}(x+2,y,z)=\varphi^{(+)}(x,y+2,z).
\end{cases}
\end{align}

The problem of determining the leading term of the asymptotic description of the front $h_0(x,y,t)$ takes the form
\begin{align}\label{eq44}
\begin{cases}
\displaystyle {h_0}_{t} =\frac{1}{2} \left( {h_0}_{x}+{h_0}_{y}-1 \right) \left( \varphi^{(-)} (x,y,h_0)+\varphi^{(+)} (x,y,h_0)\right),\\
h_0(x,y,0)=h_0^{*}=0, \quad h_0(x,y,0)=h_0(x+2,y,0)=h_0(x,y+2,0).
\end{cases}
\end{align}

\begin{figure}[H]
\begin{center}
\includegraphics[width=0.4\linewidth,height=0.4\textwidth,keepaspectratio]{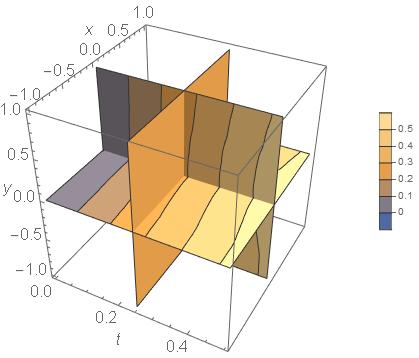}
\caption{Numerical solution of \eqref{eq44} for $t \in [0, 0.5] $.}
\label{fig:x0example1}
\end{center}
\end{figure}

From the numerical solution of \eqref{eq44}, we see that the transition layer is located in the region $-1 \leq h_0 (x,y,t) \leq 1$ for any $x\in [-1;1], y \in [-1;1]$, and $t \in [0,0.5]$ (Fig.~\ref{fig:x0example1}), thus Assumption~\ref{A3} is satisfied.

To obtain a numerical solution of \eqref{forwardexample1}, the results of which will be used for the inverse problem, it is necessary to know the initial function, which we take in the form
\begin{multline}
\displaystyle u_{init} (x,y,z,\mu )=\frac{u^{a}-u^{-a}}{2} \tanh\left(x+y+\frac{z-h_0(x,y,0)}{0.1\mu}\right)+ \frac{u^{a}+u^{-a}}{2} \\
=3\tanh\left(x+y+\frac{z}{0.001}\right)-1,
\end{multline}
with an inner transition layer in the vicinity of $h_0(x,y,0)=0$.

All the assumptions formulated in this paper are satisfied, and the considered equation \eqref{forwardexample1} has a solution in the form of an autowave with a transitional moving layer localized near $h_0(x,y,t)$, which in the zeroth-order approximation has the form
\begin{align} \label{asymptoticsolEXAMPLE1}
U_{0}=\begin{cases}
\displaystyle \varphi^{(-)}(x,y,z) +\frac{\varphi^{(+)}(x,y,h_0)-\varphi^{(-)}(x,y,h_0)}{\exp \left( \frac{ \left(h_0-z \right) \left(\varphi^{(+)}(x,y,h_0)-\varphi^{(-)}(x,y,h_0) \right) \left(1-{h_0}_{x}-{h_0}_{y} \right)}{2 \mu} \right)+1} , z \in [-1;h_0],\\
\displaystyle \varphi^{(+)}(x,y,z)+\frac{\varphi^{(-)}(x,y,h_0)-\varphi^{(+)}(x,y,h_0)}{\exp \left( \frac{ \left(h_0-z \right) \left(\varphi^{(-)}(x,y,h_0)-\varphi^{(+)}(x,y,h_0) \right) \left(1-{h_0}_{x}-{h_0}_{y} \right)}{2 \mu} \right)+1} , z \in [h_0;1].
\end{cases}
\end{align}

\begin{figure}[H]
\centering
\subfigure[]{
\includegraphics[width=0.4\linewidth]{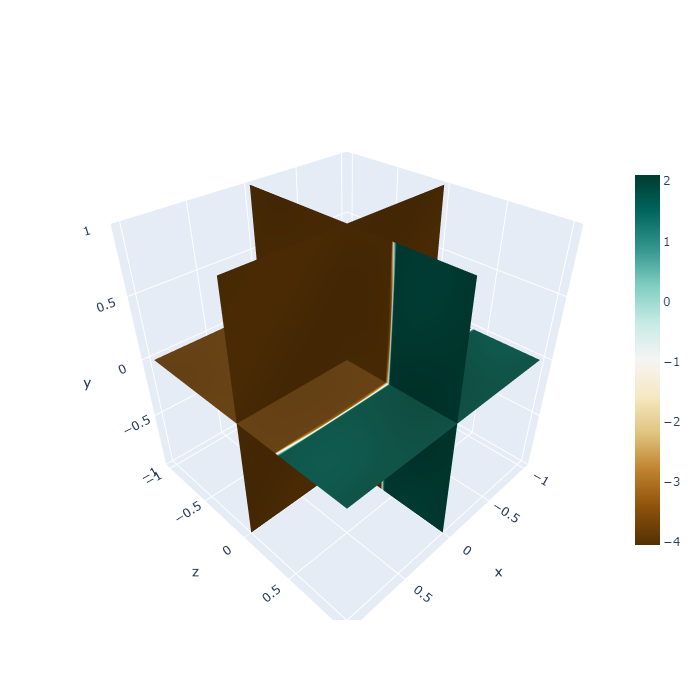} \label{fig:asymptSOL} }
\hspace{0ex}
\subfigure[]{
\includegraphics[width=0.4\linewidth]{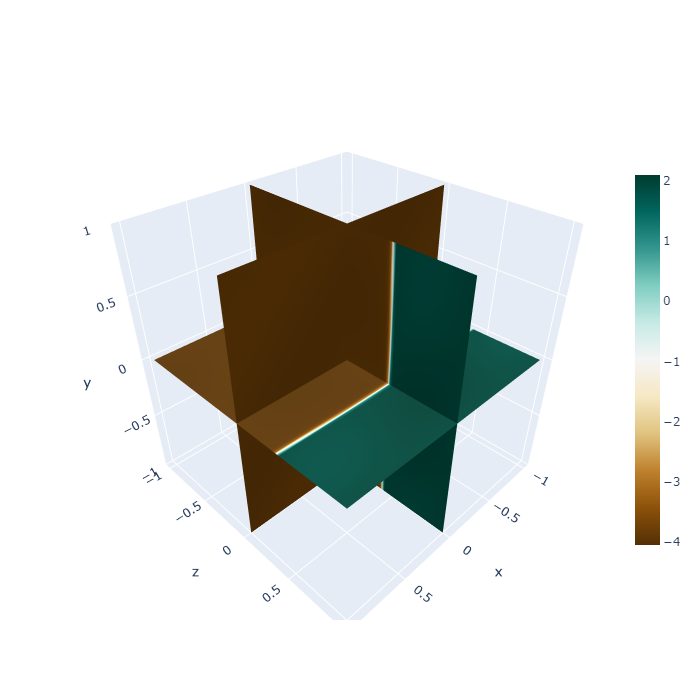} \label{fig:numSOL} }
\caption{Asymptotic solution \subref{fig:asymptSOL} and numerical solution \subref{fig:numSOL} of PDE \eqref{forwardexample1} with $t_0=0.4$. } \label{solutionsexample1}
\end{figure}

To evaluate the accuracy of our approach, we show the asymptotic solution \eqref{asymptoticsolEXAMPLE1} and a high-accuracy numerical solution\footnote{In this paper, the numerical solution is produced by a finite-volume method and can be identified as the exact solution.} at time $t_0=0.4$ in Figures~\ref{fig:asymptSOL} and \ref{fig:numSOL}, respectively, from which we conclude that the asymptotic solution is highly qualified. Indeed, the error of the asymptotic solution is
$$ \frac{\| U_0(x,y,z,t_0) - u(x,y,z,t_0) \|_{L^{2}([-1,1] \times [-1,1]\times [-1,1])}}{\|u(x,y,z,t_0) \|_{L^{2}([-1,1] \times [-1,1]\times [-1,1])}} = 0.0735. $$

\subsection{Inverse problem}

Now, as a more complex scenario, we consider the problem of identifying the source function $f(x,y,z)$ of problem \eqref{forwardexample1} when only the noisy measurements $\{u^\delta_{i,j,k}\}^{p,q,v}_{i,j,k=0}$ are provided.

Because the functions $\varphi^{(-)}(x,y,z)$ and $\varphi^{(+)}(x,y,z)$ are solutions to problems \eqref{eqUL} and \eqref{eqUR}, the asymptotic analysis allows us to reduce the original inverse problem of determining the source function $f(x,y,z)$ to a first-order equation of the asymptotic approximation of the solution of the direct problem. We use the parameters $a=1, t_0=0.4, \delta = 5\%, p=100, q=100, v=100, v^{(-)}=58 , v^{(+)}=73$.

\begin{figure}[H]
\centering
\subfigure[]{
\includegraphics[width=0.44\linewidth]{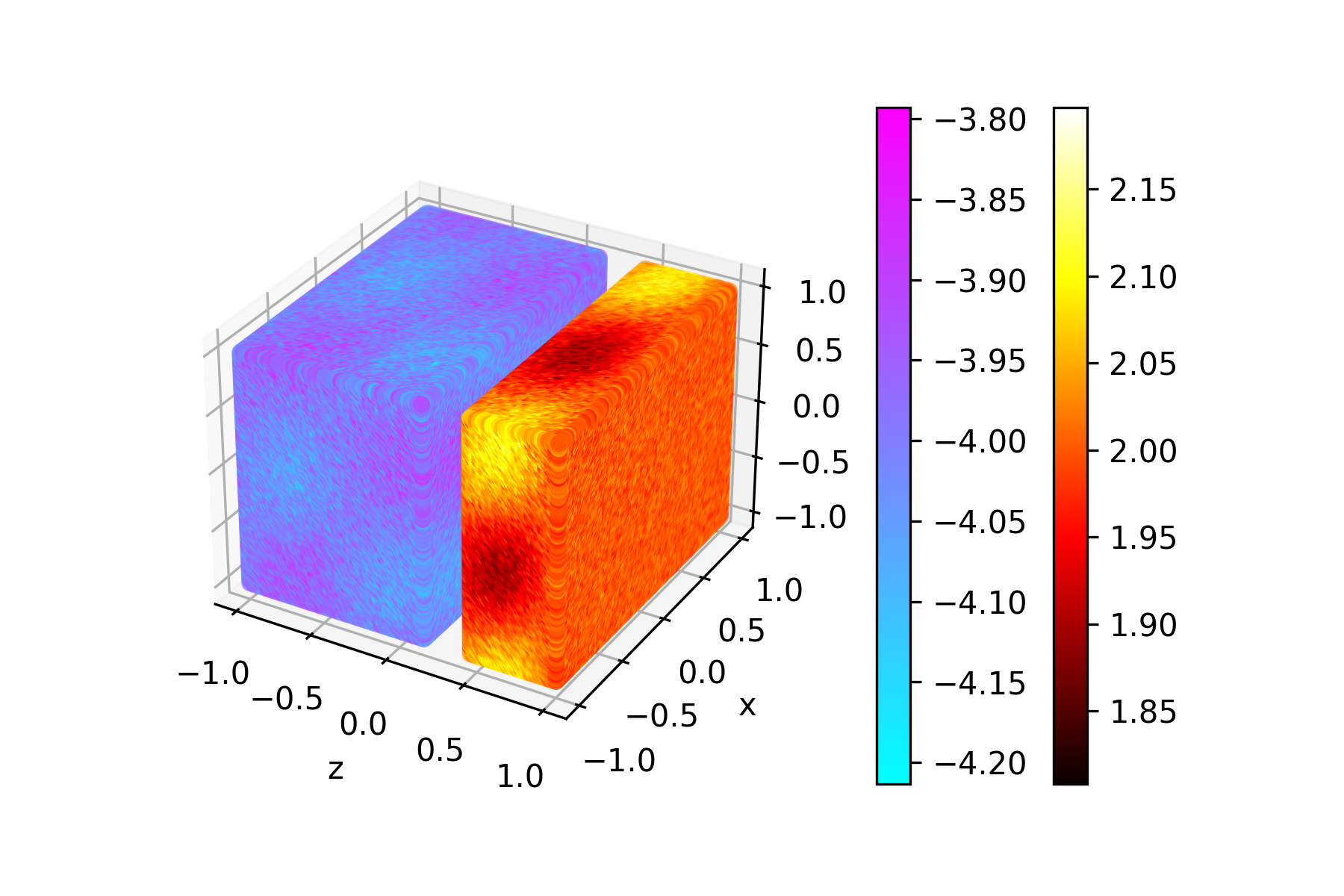} \label{fig:origNoisy} }
\hspace{0ex}
\subfigure[]{
\includegraphics[width=0.44\linewidth]{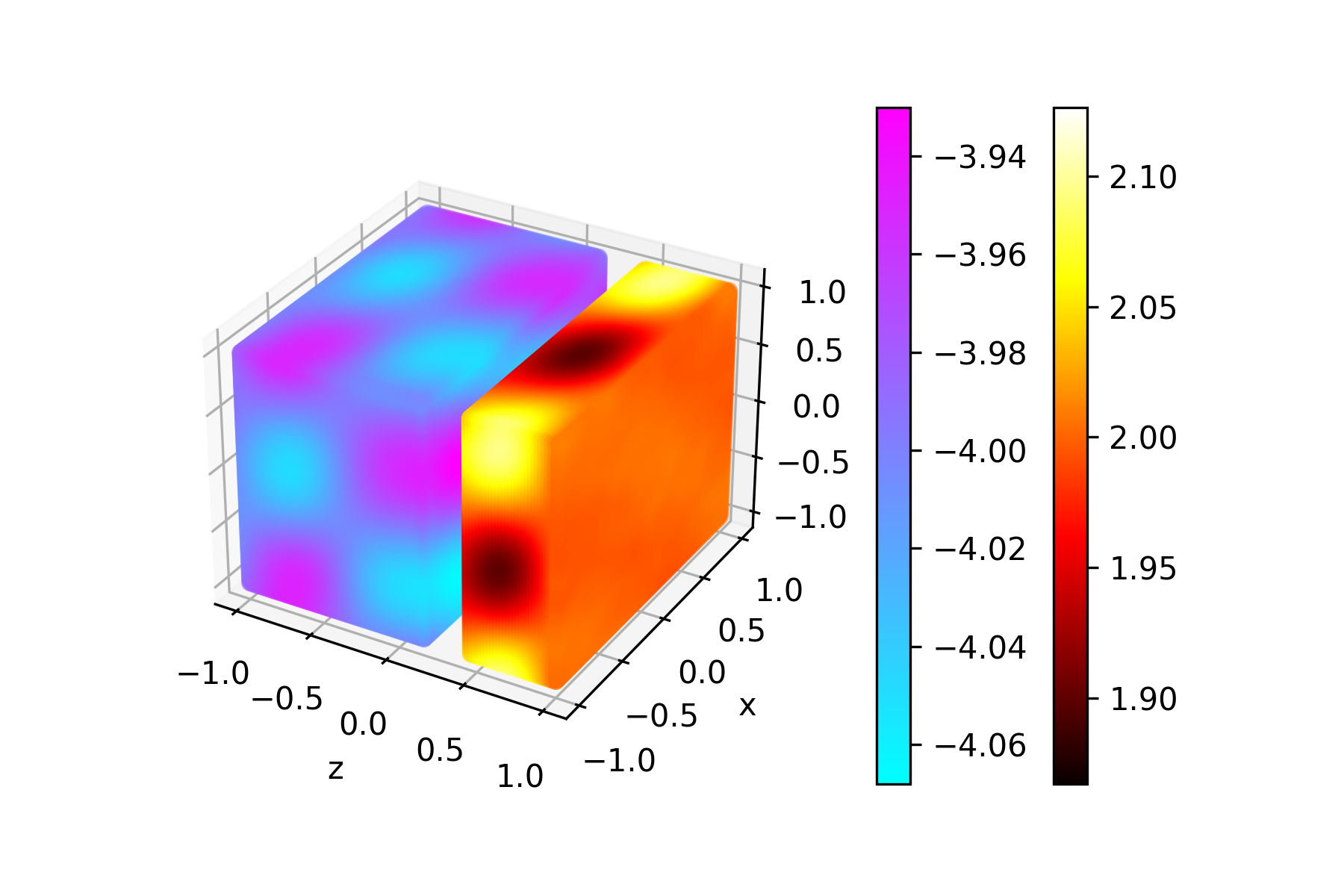} \label{fig:smoothed} }
\caption{Noisy data $\{u^{\delta}_{i,j,k}(t_0)\}$ with $\delta=5\% $ \subref{fig:origNoisy} and smooth approximate data $\{u^{\varepsilon}_{i,j,k}(t_0)\}$ \subref{fig:smoothed}. }
\label{fig:ualpha10-3}
\end{figure}

We take the data $u(x_i,y_j,z_k,t_0)$ obtained from the numerical result (i.e., by a finite-volume method) for the forward problem [Fig.~\ref{fig:numSOL}], which belong to the mesh knots $X_p=\lbrace x_i, 0 \leq i \leq p: x_i= x_0+d_1 i, d_1 = 2a/p \rbrace $, $Y_q=\lbrace y_j, 0 \leq j \leq q: y_j= y_0+d_2 j, d_2 = 2a/q \rbrace $, and $Z_v=\lbrace z_k, 0 \leq k \leq v: z_k= -a+d_3 k, d_3 = 2a/v \rbrace $. The points belonging to the transitional layer $(h_0(x,y,t_0)-\Delta h /2, h_0(x,y,t_0)+\Delta h/2)$ have been removed, so we use only the points with indices $i= 0, \cdots, p$, $j= 0, \cdots, q$ and $k= 0, \cdots, v^{(-)}$, $k=v^{(+)}, \cdots, v$. Artificial noise data are obtained by adding uniform noise with $\delta=5\%$ to the values $u(x_i, y_j, z_k, t_0)$ [Fig.~\ref{fig:origNoisy}]:
\begin{equation}
\label{noiseUni}
u^{\delta}_{i,j,k} := [1+\delta(2\, \text{rand} -1)] u(x_i,y_j,y_k,t_0),
\end{equation}
where ``rand'' returns a pseudo-random value drawn from a uniform distribution on $[0,1]$.

\begin{figure}[H]
\centering
\subfigure[]{
\includegraphics[width=0.4\linewidth]{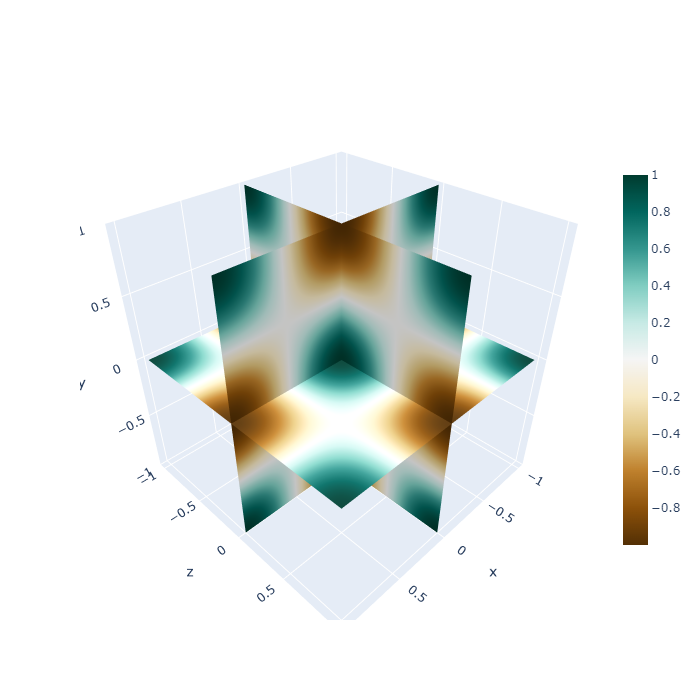} \label{fig:restSF} }
\hspace{0ex}
\subfigure[]{
\includegraphics[width=0.4\linewidth]{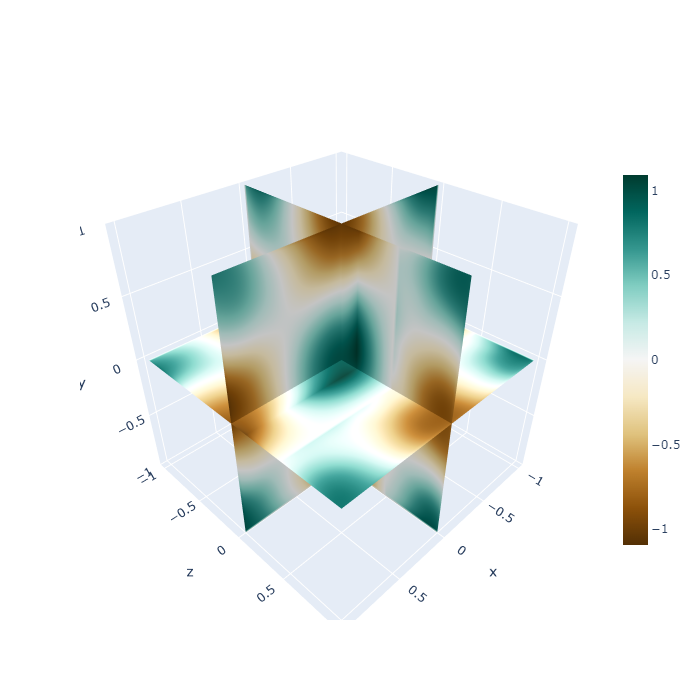} \label{fig:origSF} }
\caption{Reconstructed source function $f^{\delta}(x,y,z)$ \subref{fig:origSF} for $\delta=5\%$, where data in the transition layer was obtained by the linear interpolation, and exact source function $f^{*}(x)$ \subref{fig:restSF}. } \label{fig:sourcereconstruction}
\end{figure}

Following Algorithm~1, we smooth the noisy data on both regions by \eqref{uAlphaL} and \eqref{uAlphaR} and obtain the smooth function $u^\varepsilon(x,y,z,t_0)$ [Fig.~\ref{fig:smoothed}]. We then calculate the regularized approximate source function $f^{\delta}$ by \eqref{fdelta} and compare it with the exact source function $f^{*}(x,y,z)=\cos{ ( \pi x )} \cos{( \pi y )} \cos{( \pi z )}$ (Fig.~\ref{fig:sourcereconstruction}). The relative error of reconstructing the source function for the noise $\delta=5\%$ is $$ \frac{\| f^{\delta} - f^* \|_{L^2([-1,1] \times [-1,1]\times [-1,1])}}{\|f^* \|_{L^2([-1,1] \times [-1,1]\times [-1,1])}} = 0.3674.$$

\section{Conclusion} \label{Conclusion}

In this paper, by using asymptotic expansions, we developed an efficient approach for obtaining the approximate solutions of both forward and inverse problems associated with a three-dimensional nonlinear PDE of the RDA type with an internal transition layer. The merit of our approach is twofold. First, by incorporating the asymptotic method of inequalities, the existence and uniqueness of a classical solution of the PDE have been proven. Second, the procedure of asymptotic expansion can be used as one of model reduction for some inverse problems: when constructing an asymptotic solution, a equation linking measurements and unknown quantities can be explored, which produces an efficient inversion algorithm.

Note that our approach could be applied easily to many other similar problems in mathematical physics. For instance, the present model \eqref{mainproblem} can be viewed as a special case of the general model
$$
\mu \Delta u-\frac{\partial u}{\partial t}=\left(\mathbf{A}\left(u, x, y,z\right), \nabla\right) u+F(u, x, y,z),
$$
where $\mathbf{A}(u, x, y,z)=\left\{A_{1}(u, x, y,z), A_{2}(u, x, y,z), A_{3}(u, x, y,z)\right\}$. Having chosen the functions $A_{1}$, $A_{2}$, $A_{3}$, and $F$ so that it is possible to obtain explicit solutions for the inner and outer functions, all the results in this paper can be derived similarly.

\section*{Funding}
This work was funded by the National Natural Science Foundation of China (Grant No.\ 12171036) and the Beijing Natural Science Foundation (Key Project No.\ Z210001).

\bibliography{ChaikovskiiZhang}
\bibliographystyle{unsrt}

\end{document}

%% file: ex_shared.tex
\usepackage{lipsum}
\usepackage{amsfonts}
\usepackage{graphicx}
\usepackage{epstopdf}
\ifpdf
  \DeclareGraphicsExtensions{.eps,.pdf,.png,.jpg}
\else
  \DeclareGraphicsExtensions{.eps}
\fi

\usepackage{subfigure}
\usepackage{color}

\usepackage{longtable}
\usepackage{accents}

\usepackage{hyperref}
\usepackage{float}
\usepackage{epsfig}
\usepackage{epstopdf}
\usepackage{array,booktabs}
\usepackage{tabularx}
\usepackage{breqn}
 
\usepackage{algorithm}
\usepackage{algorithmicx}
\usepackage{algpseudocode}

 \usepackage{mathtools}
\usepackage{enumitem}
 
\newcommand{\RomanNumeralCaps}[1]
    {\MakeUppercase{\romannumeral #1}}

\usepackage{enumitem}
\setlist[enumerate]{leftmargin=.5in}
\setlist[itemize]{leftmargin=.5in}

\newsiamremark{remark}{Remark}
\newsiamremark{hypothesis}{Hypothesis}
\crefname{hypothesis}{Hypothesis}{Hypotheses}
\newsiamthm{claim}{Claim}

\newsiamthm{assumption}{Assumption}

\headers{Asymptotic expansion regularization for 3D inverse source problems}{Dmitrii Chaikovskii and Ye Zhang}

\title{Solving forward and inverse problems in a non-linear 3D PDE via an asymptotic expansion based approach \thanks{ Corresponding Author: Ye Zhang.  
\funding{The work was funded by the National Natural Science Foundation of China (No. 12171036) and Beijing Natural Science Foundation (Key project No. Z210001).}}}

\author{Dmitrii Chaikovskii\thanks{Shenzhen MSU-BIT University, 518172 Shenzhen, China 
  (\email{dmitriich@smbu.edu.cn}).}
\and  Ye Zhang\thanks{School of Mathematics and Statistics, Beijing Institute of Technology, 100081 Beijing, China 
  ( \email{ye.zhang@smbu.edu.cn}).} 
}

\usepackage{amsopn}
